\documentclass[noams,lp]{compositio}

\usepackage{amsbsy,amssymb,amscd,amsfonts,latexsym,amstext,delarray,
 amsmath,stackrel,lineno,tabularx,url,cite}

\input xypic
\usepackage[pdfborder={0 0 0}]{hyperref}

\hypersetup{
    colorlinks = true,
    linkcolor = blue,
    anchorcolor = red,
    citecolor = green,
    filecolor = red,
    pagecolor = red,
    urlcolor = magenta
}

\newcommand*\patchAmsMathEnvironmentForLineno[1]{%
  \expandafter\let\csname old#1\expandafter\endcsname\csname #1\endcsname
  \expandafter\let\csname oldend#1\expandafter\endcsname\csname end#1\endcsname
  \renewenvironment{#1}%
     {\linenomath\csname old#1\endcsname}%
     {\csname oldend#1\endcsname\endlinenomath}}%
\newcommand*\patchBothAmsMathEnvironmentsForLineno[1]{%
  \patchAmsMathEnvironmentForLineno{#1}%
  \patchAmsMathEnvironmentForLineno{#1*}}%
\AtBeginDocument{%
\patchBothAmsMathEnvironmentsForLineno{equation}%
\patchBothAmsMathEnvironmentsForLineno{align}%
\patchBothAmsMathEnvironmentsForLineno{flalign}%
\patchBothAmsMathEnvironmentsForLineno{alignat}%
\patchBothAmsMathEnvironmentsForLineno{gather}%
\patchBothAmsMathEnvironmentsForLineno{multline}%
}

\newtheorem{thm}{Theorem}[section] 
\newtheorem{defn}[thm]{Definition} 
\newtheorem{prop}[thm]{Proposition}
\newtheorem{cor}[thm]{Corollary}
\newtheorem{lem}[thm]{Lemma}

\newtheorem{rem}[thm]{Remark}
\newtheorem{example}[thm]{Example}

\def\Aut{{\rm Aut}}

\def\Hom{{\rm Hom}}

\def\ssh{{\rm Sh}}

\def\N{{\mathbb N}}

\def\R{{\mathbb R}}
\def\Z{{\mathbb Z}}

\def\cC{{\mathcal C}}

\def\cI{{\mathcal I}}
\def\cF{{\mathcal F}}
\def\cJ{{\mathcal J}}

\def\cT{{\mathcal T}}

\def\qqq{\,,\quad~\forall}

\newcommand{\ie}{{\it i.e.\/}\ }

\newcommand{\cf}{{\it cf.}}
\newcommand{\opcit}{{\it op.cit.\/}\ }

\def\Hom {{\mbox{Hom}}}

\def\cF{\mathfrak F}

\def\Se{\frak{ Sets}}
\def\fin{\frak{ Fin}}
\def\yon{\frak{ y}}
\def\Arc{\frak{Arc}}
\def\tt{\frak{t}}

\begin{document}

\title{Cyclic structures and the topos of simplicial sets}
\author{Alain Connes}
\email{alain@connes.org}
\address{Coll\`ege de France,
3 rue d'Ulm, Paris F-75005 France\newline
I.H.E.S. and Ohio State University.}
\author{Caterina Consani}
\email{kc@math.jhu.edu}
\address{Department of Mathematics, The Johns Hopkins
University\newline Baltimore, MD 21218 USA.}
%
\thanks{The second author is partially supported by the NSF grant DMS 1069218
and would like to thank the Coll\`ege de France for some financial support.
}

\begin{abstract}
Given a  point $p$ of the topos $\hat\Delta$ of simplicial sets and the corresponding flat covariant functor $\cF:\Delta\longrightarrow\Se$, we determine the extensions of $\cF$   to the cyclic category $\Lambda\supset \Delta$. We show that to each such cyclic structure on a point $p$ of $\hat\Delta$ corresponds a group $G_p$, that such groups can be noncommutative and that each $G_p$ is described as the quotient of a left-ordered group by the subgroup generated by a central element. Moreover for any cyclic set $X$ the fiber (or geometric realization) of the underlying simplicial set of $X$ at $p$ inherits canonically the structure of a $G_p$-space. This gives a far reaching generalization of the well-known circle action on the geometric realization of cyclic sets.
\end{abstract}

\maketitle

\tableofcontents  

\section{Introduction}

This paper aims to illustrate the unifying power of the notion of topos, due to Grothendieck, in the context of cyclic homology. Our main motivation originates from the recent discovery (\cf\cite{cycarch}) of the role of cyclic homology of schemes for the cohomological interpretation of the archimedean local factors of L-functions of arithmetic varieties. This result raises naturally the question of a conceptual interpretation of cyclic homology of schemes.  In \cite{CoExt} it was shown that cyclic cohomology can be interpreted as Ext-functor in the category of cyclic modules. These modules are defined as contravariant functors from the (small) cyclic category $
\Lambda$ to the category of abelian groups. This development brings at the forefront the crucial role played by the cyclic category as an extension of the simplicial category $\Delta$ by finite cyclic groups. Moreover, in \opcit it was also shown   that  the classifying space $B\Lambda$ is equal to the classifying space of the topological group $U(1)$, and later  it was  discovered (\cf\cite{good, bu0, jones}) that the geometric realization of the simplicial set underlying a cyclic set (this latter understood as a contravariant functor $\Lambda\longrightarrow\Se$ to the category of sets) inherits naturally an action of  $U(1)$.   The equality $B\Lambda=BU(1)$ then  leads to a deep analogy between cyclic cohomology and $U(1)$-equivariant cohomology. 

In this article we show that one obtains a conceptual understanding of the above $U(1)$-action on the geometric realization of cyclic sets by extending that result to the framework of topos theory.  The transition from a small category to its classifying space produces in general a substantial loss of information: the classifying space of $\Delta$ is, for instance, a singleton. It is exactly at this point that the implementation of topos theory turns out to be useful  to provide the correct  environment that encloses both schemes and  small categories while also furnishing the tools for the development of cohomology. 
In the words of Grothendieck:

{\it 
L'id\'ee du {\it topos}  englobe, dans une intuition topologique commune, aussi bien
les traditionnels espaces (topologiques), incarnant le monde de la grandeur continue, que
les (soi-disant) ``espaces'' (ou ``vari\'et\'es'') des g\'eom\`etres alg\'ebristes abstraits imp\'enitents,
ainsi que d'innombrables autres types de structures, qui jusque-l\`a avaient sembl\'e riv\'ees
irr\'em\'ediablement au ``monde arithm\'etique'' des agr\'egats ``discontinus'' ou ``discrets''.}

The category $\ssh(X)$ of sheaves of sets on a topological space $X$ is a topos that captures all the relevant information
on $X$.
For a small category $\cC$, the associated category $\hat \cC=\Se^{\cC^{\rm op}}$ of contravariant functors $\cC\longrightarrow \Se$ is a topos. What is more, the notion of point is meaningful for any topos $\cT$: a point is simply a geometric morphism $f:\Se\longrightarrow \cT$ from the topos  of sets to $\cT$. To each point  of  $\cT$ corresponds a contravariant functor $\cT\longrightarrow\Se$ which is the inverse image part of the geometric morphism $f$ and that preserves finite limits and arbitrary colimits. This picture generalizes the functor that associates  to a sheaf of sets on a topological space $X$  the stalk at a point of $X$. In particular, for the topos $\hat \cC=\Se^{\cC^{\rm op}}$ ($\cC$ a small category), the points are described by  {\em flat, covariant} functors $\cF:\cC\longrightarrow\Se$. Then the inverse image part of the geometric morphism associated to a point of $\hat \cC$ determines a natural generalization of the notion of the geometric realization  of a simplicial set. This latter notion is obtained, in the case $\cC=\Delta$, by considering  the flat functor $\cF=\underline \Delta: \Delta \longrightarrow \Se$ that associates to an integer $n\geq 0$ the standard $n$-simplex. In general,  the flatness of $\cF$  implies that the geometric realization functor
\begin{equation*}
\label{gerea}
|~|:\Se^{\cC^{\rm op}}\longrightarrow \Se,\quad R \mapsto |R|:=R\otimes_\cC \cF
\end{equation*}
is left exact thus transforming finite products in $\Se^{\cC^{\rm op}}$ into finite products in $\Se$. This property combines with the fact that $\cF=\underline \Delta$ extends to the larger cyclic category $\Lambda$ to yield the natural action of the group  $U(1)$ on the geometric realization of a cyclic set. 

In this paper we provide a far reaching generalization of this construction for the points of the topos of simplicial sets.  Given a point $p$ of the topos $\hat \Delta$ with associated flat functor $\cF:\Delta\longrightarrow \Se$, a {\em cyclic structure} on $p$ is defined to be an extension of $\cF$ from $\Delta$ to $\Lambda$. Our main result states that  cyclic structures are classified by the datum provided by a totally left-ordered group $G$  (not necessarily abelian) endowed with a central element $z>1$ such that the interval $[1,z]\subset G$ generates $G$.  It is well known (\cf\cite{MM}) that the points of the topos $\hat \Delta=\Se^{\Delta^{\rm op}}$ correspond to intervals $I$ \ie  totally ordered sets with a minimal element $b$ and a maximal element $t>b$. Proposition \ref{propcp} shows
that if $(G,z)$ is a left-ordered group with a fixed central element $z>1$, one obtains a natural cyclic structure on the point $p_I$ of $ \hat\Delta$ associated to the interval $I=[1,z]$.  The converse of this statement constitutes the main result of this paper. More precisely one shows (\cf~ Theorem \ref{main})  the following 

\textbf{Theorem}
~{\em Let $p$  be a point of the topos of simplicial sets and let $I$ be the associated interval. Let
$ G=\left(\Z\times I\right)/\sim$~  be endowed with the lexicographic ordering, where the equivalence relation identifies  $(n,b)\sim (n-1,t)$, $\forall n\in \Z$. Then, a cyclic structure on $p$ corresponds to a group law on $G$ such that\newline
$(a)$~ the order relation on $G$ is left invariant \newline
$(b)$~the following equalities hold
   \[
  (n,b)(m,u)=(m,u) (n,b)=(n+m,u),\qquad \forall n,m\in \Z, ~u\in I.
   \]}
The proof of the above theorem  is rather involved and  is developed in \S \ref{mainsect}.
As an immediate corollary one obtains examples (Corollaries \ref{hom} and \ref{hom1}) of points $p = p_I$ of the topos $\hat \Delta$ with a prescribed number of cyclic structures. 
Moreover, Theorem \ref{thmgeom1}
shows that under the hypothesis of the above theorem, there  exists, for any cyclic set $S$, a {\em canonical} right action of the quotient group $G/\Z$, of $G$ by the central subgroup $\{(n,b)\mid n\in \Z\}\sim \Z$, on 
the ``geometric realization"  $|S|_p=p^*S$ corresponding to the point $p$ of the topos $\hat\Delta$. \newline
In the last part of the paper (\cf~\S\ref{the last}) we discuss the subtlety that distinguishes the notion of cyclic structure on a point of $\hat\Delta$ as given here, from the definition of a point of the topos of cyclic sets (called ``abstract circle") as in the unpublished note \cite{Moerdjik}. By using \opcit we prove  that the category of points of the topos of cyclic sets is {\em equivalent} to the category $\Arc$ of archimedean sets whose objects are pairs $(X,\theta)$ made by a totally ordered set $X$ together with an automorphism  $\theta$ of $X$ such that $\theta(x)>x$, $\forall x\in X$ and fulfilling the following archimedean property: for any given pair $x,y\in X$ there exists $n\in \N$ such that $y<\theta^n(x)$. Modulo the identification $f = \theta^m\circ f$ ($f$ morphism in $\Arc$: see  Definition~\ref{last}), one obtains the equivalence of $\Arc$ with the category of abstract circles: \cf~Proposition~\ref{equivcat}. To the inclusion of categories $\Delta\subset \Lambda$ corresponds a geometric morphism of topoi and a related map between the corresponding points that associates to an interval $I$ the archimedean set $X=(\Z\times I)/\sim$ as in  Theorem \ref{main}, endowed with the translation
$\theta(n,u)=(n+1,u)$. In this way, one derives a reformulation of Theorem \ref{main} by stating that the cyclic structures on a point of $\hat\Delta$ are classified by  (order compatible) group structures on the associated point of the topos of cyclic sets. \newline
Thus the notion of cyclic structure on a point of the topos of simplicial sets is a very subtle concept involving both groups and linear orders, not  be confused with the definition of point of the topos of cyclic sets, that only involves linear order.

\section{The simplicial and the cyclic categories. The cyclic structures}

In this section we recall the definitions of the simplicial and the cyclic categories together with the definition of points of the related topoi. Then, in \S\ref{cycpoint} we introduce the new notion of cyclic structure on the points of the topos of simplicial sets.

\subsection{The  simplicial category $\Delta$}

We recall the classical presentation of the simplicial category $\Delta$.
Let $\Delta$ be the small category with one object $[n]=\{0,\ldots ,n\}$ for each non-negative integer $n\geq 0$: the object $[n]$ is viewed as a totally ordered set.  The morphisms in $\Delta$ are non decreasing maps of sets
\begin{equation*}
    \Hom_\Delta([n],[m])=\{f:\{0,\ldots ,n\}\to \{0,\ldots ,m\}\mid x\geq y\implies f(x)\geq f(y)\}.
\end{equation*}
The following well-known result (\cf\cite{NCG} III.A.$\alpha$, Proposition 2) implies a presentation of $\Delta$ by generators and relations
\begin{prop}\label{presdelta}
$(i)$~For $j\in \{0,\ldots ,n\}$, there are $n+1$ surjective morphisms
\begin{equation*}
\sigma_j\in \Hom_\Delta([n+1],[n]),\quad \sigma_j(i)=\left\{
                                   \begin{array}{ll}
                                     i& \hbox{if}~\ i\leq j\\
                                     i-1 & \hbox{if}~\ i>j.
                                   \end{array}
                                 \right.
\end{equation*}
$(ii)$~For $j\in \{0,\ldots ,n\}$, there are $n+1$ injective morphisms
\begin{equation*}
\delta_j\in \Hom_\Delta([n-1],[n]),\quad    \delta_j(i)=\left\{
                                   \begin{array}{ll}
                                     i& \hbox{if}~\ i< j\\
                                     i+1 & \hbox{if}~\ i\geq j.
                                   \end{array}
                                 \right.
\end{equation*}
$(iii)$~The following relations hold in $\Delta$
\begin{equation}\label{relj}
  \sigma_j\circ \sigma_i=\sigma_i\circ\sigma_{j+1}\qquad {\rm for}\quad i\leq j, \qquad  \delta_j\circ \delta_i=\delta_i\circ \delta_{j-1}\qquad {\rm for}\quad i<j,
\end{equation}
\begin{equation*}
  \sigma_j\circ \delta_i=\left\{
                           \begin{array}{ll}
                            \delta_i\circ \sigma_{j-1} & \hbox{if}~\ i<j \\
                             {\rm id} & \hbox{if}~\ i=j\ \hbox{or}\ i=j+1 \\
                            \delta_{i-1}\circ \sigma_j & \hbox{if}~\ i>j+1.
                           \end{array}
                         \right.
\end{equation*}
\end{prop}
Any morphism in $\Delta$ can be written as  a composite $(\prod_j \delta_j\circ\prod_k \sigma_k)$ of a (finite) product of $\delta$'s with a (finite) product of $\sigma$'s. Moreover, by implementing the relations \eqref{relj}, one can also re-order these products  so that the indices for the $\sigma$'s are increasing  while they are non-decreasing for the $\delta$'s.

\subsection{The cyclic category $\Lambda$}\label{cyccat}

The cyclic category $\Lambda$ was introduced in \cite{CoExt} using classes of degree-one maps of the circle. In this paper we use a reformulation of $\Lambda$ that implements the lift of these maps to the universal cover (\cf~\cite{FT}). The advantage of this description is to involve only {\em periodic non-decreasing} maps $\Z\to \Z$.
For each pair of integers $a>0$, $b>0$, consider the set
\begin{equation*}
    \cI(a,b)=\{f:\Z\to \Z\mid f(x)\geq f(y),~\forall x\geq y; \ \ f(x+a)=f(x)+b,~ \forall x\in \Z\}.
\end{equation*}
The following relation is an equivalence relation on $\cI(a,b)$
\begin{equation}\label{equiv}
   f\sim g\iff \exists k\in \Z, \ g(x)=f(x)+kb\qquad \forall x\in \Z.
\end{equation}
Denote by $\cC(a,b)$  the set of the equivalence classes of maps in $\cI(a,b)$.
\begin{prop}
$(i)$~The set $\cC(a,b)$ is finite.

$(ii)$~The equivalence relation \eqref{equiv} is compatible with the composition of maps.
\end{prop}
\proof $(i)$~In any given equivalence class in $\cC(a,b)$ there is a unique map $f: \Z \to \Z$ such that $f(0)\in \{0,\ldots, b-1\}$, since $\{0,\ldots, b-1\}$ is a fundamental domain for the subgroup $b\Z\subset\Z$. It follows that $f(a)=f(0)+b\in \{b,\ldots, 2b-1\}$, thus there are only finitely many possibilities for $f(x)\in [f(0),f(a)]$ as $x\in \{0,\ldots, a-1\}$.  One concludes, by applying the periodicity property, that there can be only finitely many  equivalence classes in each $\cC(a,b)$.

$(ii)$~Let $j\in\{1,2\}$ be an index and let $f_j\in\cI(a,b)$ and $g_j\in \cI(b,c)$. One has $g_1\sim g_2~\Leftrightarrow~g_2(x)=g_1(x)+kc$, for some $k\in\Z$ and $\forall x\in \Z$. Then it follows that $g_2\circ f_j\sim g_1\circ f_j$.
Similarly one has $f_1\sim f_2~\Leftrightarrow~f_2(x)=f_1(x)+k_1b$, for some $k_1\in\Z$ and $\forall x\in \Z$. It follows that $g_j\circ f_2(x)=g_j(f_1(x)+k_1b)=g_j(f_1(x))+k_1c \qqq x\in \Z$, 
thus $g_j\circ f_2\sim g_j\circ f_1$.\endproof

\begin{defn}\label{defncatcy}
The {\em cyclic category} $\Lambda$ has one object $[n]$ for each non-negative integer $n\geq 0$. The morphisms in $\Lambda$ are given by
\begin{equation*}
    \Hom_\Lambda([n],[m])=\cC(n+1,m+1).
\end{equation*}
\end{defn}
By construction one derives a natural functor (not faithful)
\[
\mu: \Lambda\longrightarrow\fin
\]
to the category $\fin$ of finite sets. $\mu$ is the identity on the objects of $\Lambda$ and  to a morphism $f\in \Hom_\Lambda([n],[m])$ it associates the map of finite sets $\mu(f): \{0, \ldots, n\} \to\{0,\ldots, m\}$ given by
\begin{equation}\label{under}
    \mu(f)(x)=f(x) \ \text{mod.} \ m+1\qqq x\in \Z/(n+1)\Z.
\end{equation}
There is an inclusion functor $j:\Delta\hookrightarrow \Lambda$ that is the identity on the objects of the simplicial category $\Delta$ and on the morphisms $f\in \Hom_\Delta([n],[m])$ is defined as follows
\begin{equation}\label{embedj}
    j(f)(x+k(n+1))= f(x)+k(m+1)\qqq x\in[n],~ k\in\Z.
\end{equation}
One easily checks that $j(f)^{-1}([0,m])=[0,n]$, since by construction one has  $j(f)([0,n]+k(n+1))\subset [0,m]+k(m+1)$, $\forall k\in \Z$. A similar result holds for all morphisms in $\Lambda$ as the following proposition (\cf part (i)) shows

\begin{prop}\label{decolam} Let $f\in \Hom_\Lambda([n],[m])$, then
\begin{enumerate}
\item For any interval $I=[y,y+m]$, $f^{-1}(I)$ is an interval of the form $[x,x+n]$.
\item There is a unique decomposition of the form
\begin{equation*}
    f=j(h)\circ t, \quad  h\in \Hom_\Delta([n],[m]), \ t\in \Aut_\Lambda([n]).
\end{equation*}
\end{enumerate}
\end{prop}
\proof $(i)$~$f^{-1}(I)$ is a finite interval $J$ because $f$ is  increasing and $f(x)\to \pm \infty$ when $x\to \pm\infty$. Moreover  since for any fixed interval $I$ as in (i) and for $k\in \Z$ (varying) the translates $I+k(m+1)$  form a partition of $\Z$, the same statement holds for their inverse images which are, in view of the periodicity property of the maps in $\Lambda$, of the form $J+k(n+1)$.\newline
$(ii)$~It follows from the last sentence in the proof of $(i)$ that there exists a unique $t\in \Aut_\Lambda([n])$ such that $f^{-1}([0,m])=t^{-1}([0,n])$. Thus $(f\circ t^{-1})^{-1}([0,m])=[0,n]$ and by implementing the functor $j$ as in \eqref{embedj} one concludes that  there exists $h\in \Hom_\Delta([n],[m])$ such that $f\circ t^{-1}=j(h)$. \endproof
One derives the following presentation of the cyclic category $\Lambda$ that in applications is usually referred to as the decomposition $\Lambda = \Delta {\bf C}$, where ${\bf C}$ is the category with the same objects as $\Lambda$  and whose morphisms are automorphisms (\cf\cite{CoExt}).

\begin{prop}\label{prescyclc}
The cyclic category $\Lambda$ admits the following presentation as an extension of the small category $\Delta$ by means of a new generator $\tau_n\in C_{n+1}:=\Aut_\Lambda([n])$ for each $n\geq 0$ and the relations
\begin{align}\label{pres1}
&&\tau_n^{n+1}=id,&\notag\\
\tau_n\circ \sigma_0&=\sigma_n\circ \tau_{n+1}^2 &&&   \tau_n\circ \sigma_j &=\sigma_{j-1}\circ \tau_{n+1} \qqq j\in \{1, \ldots, n\}\\\label{pres2}
\tau_n\circ \delta_0 &=\delta_n  &&&  \tau_n\circ \delta_j&=\delta_{j-1}\circ \tau_{n-1} \qqq j\in \{1, \ldots, n\}
\end{align}
\end{prop}
\proof We introduce the cyclic permutation $\tau$:
\begin{equation*}
\tau: \Z \to \Z\qquad    \tau(x)=x-1\qqq x\in \Z.
\end{equation*}
Note that, for each $n\geq 0$, $\tau$ yields an element $\tau_n\in C_{n+1}=\Aut_\Lambda([n])$ (the index $n$ will be dropped when the context is clear). The relation $\tau_n^{n+1}=id$ follows from the definition of the equivalence relation \eqref{equiv}. To prove \eqref{pres1} and \eqref{pres2} one first determines, using Proposition \ref{decolam} and for
$f\in \Hom_\Delta([m],[n])$, an integer $a$ such that
\begin{equation*}
(\tau_n\circ j(f))^{-1}([0,n])=\tau^{a}([0,m]).
\end{equation*}
Then it follows that $\tau_n\circ j(f)\circ \tau^a\in j(\Delta)$, thus there exists $h\in \Hom_\Delta([m],[n])$ such that $\tau_n\circ j(f)= j(h)\circ \tau^{-a}$. Since the powers of $\tau$ are automorphisms, the surjectivity {\em resp.} the injectivity properties of $f$ are inherited from those of the associated $h$. Proposition \ref{presdelta} implies that the commutation relations
\eqref{pres1} and \eqref{pres2} stay within the $\sigma$'s and $\delta$'s. Then one easily checks directly these relations. Conversely, these relations allow one to write any morphism in $\Lambda$ as a product $(j(f)\circ t)$ for some automorphism $t$ and this suffices to present the small category $\Lambda$. \endproof
We refer to \cite{Loday} ({\em cf.} \S 6.1) for a detailed proof of the above proposition and for the general notion of crossed simplicial group that is beneath and generalizes the above decomposition $\Lambda = \Delta {\bf C}$. 

\subsection{The points of the topos $\Se^{\cC^{\rm op}}$}

Let $\cC$ be a small category. We recall  that a {\em point} of the topos  $\hat \cC=\Se^{\cC^{\rm op}}$ is completely determined by a covariant functor $\cF:\cC\longrightarrow\Se$ that has the further property to be {\em flat} (\cf\cite{MM} VII.5, Definition 1 and Theorem 2). More generally, a point of a topos $\cT$ is defined as a geometric morphism $f: \Se\longrightarrow\cT$. To such  $f$ corresponds naturally the pullback functor $f^*: \cT \longrightarrow \Se$ which is the inverse image part of the geometric morphism (\opcit VII.1, Definition 1). Besides being a left adjoint functor (to $f_*$, the direct image part of $f$) and hence (right exact and) commuting with arbitrary direct limits (\ie colimits), the pull-back $f^*$ is also left exact thus it commutes with finite inverse limits. In the specific case of the topos $\cT=\hat\cC$ associated to a small category $\cC$, $f^*: \hat\cC\longrightarrow \Se$ pulls back contravariant functors $F:\cC\longrightarrow\Se$, hence the inverse image $f^*F$  is a set. In this way $f^*$ is understood as a covariant functor
 \begin{equation*}
    f^*:\Se^{\cC^{\rm op}}\longrightarrow \Se.
 \end{equation*}
 To describe the covariant functor $\cF: \cC \longrightarrow \Se$ associated to a point $f: \Se\longrightarrow \cC$ of $\hat\cC$, one uses the Yoneda embedding
\begin{equation}\label{yoneda}
   \yon: \cC\longrightarrow \Se^{\cC^{\rm op}}\qquad C\mapsto \Hom_\cC(\cdot\, , C)
\end{equation}
and defines $\cF$ as the composite:
\begin{equation}\label{yoneda1}
    \cF=f^*\circ \yon.
\end{equation}
The obtained $\cF: \cC \longrightarrow \Se$ is  a {\em filtering functor} (\cite{MM} VII.6, Definition 2). This means that the category $\cI=\int_\cC \cF$,  whose objects are pairs $(x, C)$ where $C$ is an object of $\cC$ and $x\in \cF(C)$, is a {\em filtering category} \ie $\cI$ fulfills the following properties
\begin{enumerate}
  \item $\cI$ is non empty; \ie, $\cF(C)\neq\emptyset$ for at least one object $C$ of $\cC$.
  \item For any two objects $i,j$ in $\cI$ there is a diagram $i\leftarrow k\rightarrow j$, for some object $k$ in $\cI$.
  \item For any two arrows $i\rightrightarrows j$ in $\cI$  there exists an object $k$ in $\cI$ and a commutative diagram of the form  $k\to i\rightrightarrows j$.
\end{enumerate}
Conversely, one can show that given a filtering functor $\cF:\cC\longrightarrow\Se$, the associated point $f: \Se\longrightarrow \cC$ of $\hat\cC$ is  described by the pair of adjoint functors
\begin{equation}\label{tens}
    f^*(R)=R\otimes_\cC \cF,\quad f_*(S)=\Hom_\cC(\cF,S).
\end{equation}
We refer to VII.6, Theorem 3  of \cite{MM} for a complete proof.

\subsection{The points of the topos $\hat \Delta$}

The points of the topos $\hat \Delta$ of simplicial sets are constructed from {\em intervals}, \ie totally ordered sets $I$ with a smallest element $b_I\in I$ ($b_I\leq a$, $\forall a\in I$) and a largest element $t_I\in I$: \cite{MM} VIII.8.
The intervals $I$ and the morphisms
\begin{equation*}
   \Hom_\geq(I,J)=\{f:I\to J\mid x\leq y\implies f(x)\leq f(y), \ f(b_I)=b_J, \ f(t_I)=t_J\}
\end{equation*}
given by non-decreasing maps preserving the two end points form a category $\cJ$.
One identifies the opposite category $\Delta^{\rm op}$ with the full subcategory of $\cJ$ defined by the intervals of the form $
n^*:=\{0,1,\ldots,n+1\}$, for $n\geq 0$. 
The morphisms $s_j=\sigma_j^*$ and $d_j=\delta_j^*$ in the full subcategory are given, for $j\in \{0,\ldots,n\}$, by the following formulas
\begin{equation*}
    s_j:n^*\to (n+1)^*\qquad \ s_j(i)=\left\{
                                   \begin{array}{ll}
                                     i & \hbox{if}\ i\leq j\\
                                     i+1 & \hbox{if}\ i>j
                                   \end{array}
                                 \right.
\end{equation*}
\begin{equation*}
    d_j:n^*\to (n-1)^*\qquad \ d_j(i)=\left\{
                                   \begin{array}{ll}
                                     i & \hbox{if}\ i\leq j\\
                                     i-1 & \hbox{if}\ i>j.
                                   \end{array}
                                 \right.
\end{equation*}
The point $p_I$ of $\hat \Delta$ associated to an interval $I$ is defined as in \eqref{yoneda1}  by a covariant functor $\cF_I:\Delta\longrightarrow \Se$. This functor is described as the following contravariant functor
\begin{equation}\label{fsubi}
\cF_I: \Delta^{\rm op}\longrightarrow\Se,\qquad    \cF_I(n^*):=\Hom_\geq(n^*,I).
\end{equation}
Notice that an element $\beta\in \cF_I(n^*)$ is encoded by an increasing  sequence $(\beta_j)_{0\leq j\leq n+1}$, $\beta_j\in I$, $\beta_0=b_I$, $\beta_{n+1}=t_I$. The covariant action of $\Delta$ is expressed by the following maps
\begin{equation*}
    \sigma_j:\cF_I((n+1)^*)\to \cF_I(n^*)\qquad  \sigma_j(\beta)_i=\left\{
                                   \begin{array}{ll}
                                     \beta_i & \hbox{if}\ i\leq j\\
                                     \beta_{i+1} & \hbox{if}\ i>j
                                   \end{array}
                                 \right.
\end{equation*}
 and
\begin{equation*}
    \delta_j:\cF_I((n-1)^*)\to \cF_I(n^*)\qquad\delta_j(\beta)_i=\left\{
                                   \begin{array}{ll}
                                     \beta_i & \hbox{if}\ i\leq j\\
                                     \beta_{i-1} & \hbox{if}\ i>j.
                                   \end{array}
                                 \right.
\end{equation*}
Notice that $\sigma_j$ removes $\beta_{j+1}$ from the list whereas $\delta_j$ repeats $\beta_j$.

\subsection{The cyclic structures on points of $\hat \Delta$}\label{cycpoint}

We are now ready to introduce the central notion of cyclic structure that will be further developed in the remaining sections of this paper.

\begin{defn}\label{defncycpoints}
Let $p$ be a point of the topos $\hat \Delta$ of simplicial sets and let $\cF=p^*\circ \yon: \Delta \longrightarrow \Se$ be the associated (filtering) functor. We call a cyclic structure on $p$ the datum provided by an extension of $\cF$ to a functor $\tilde \cF:\Lambda \longrightarrow \Se$.
\end{defn}

 Note that the same point of $\hat\Delta$ can have different cyclic structures and that there are points without any cyclic structure: we refer to \S\ref{sectcycpoints}, Corollaries \ref{hom} and \ref{hom1}.

\begin{example}
{\rm  As a basic example of cyclic structure we consider the point $p_{[0,1]}$ of $\hat \Delta$ associated to the interval $[0,1]\subset \R$. First we describe the functor $\cF_{[0,1]}$ as in \eqref{fsubi} in an equivalent way by encoding the datum of the increasing sequence $(\beta_j)_{0\leq j\leq n+1}$, $\beta_j\in [0,1]$, $\beta_0=0$, $\beta_{n+1}=1$, by means of the point of the standard $n$-simplex 
\begin{equation}\label{simplex}
    \sigma(\beta)=\sum_{j=0}^n (\beta_{j+1}-\beta_j)v_j\in \underline \Delta^n
\end{equation}
where the $v_j$'s are linearly independent chosen vectors of a fixed, real vector space $V$ and $\underline \Delta^n$ denotes the standard $n$-simplex in $V$. In this way we obtain an equivalence of the functor $\cF_{[0,1]}$ with the covariant functor  $\underline \Delta:\Delta\longrightarrow\Se$ that associates to the object $[n]$ of $\Delta$  the $n$-simplex $\underline \Delta^n$ and to the morphism $f\in \Hom_\Delta([n],[m])$ its {\em unique} extension as an affine map $\tilde f: V\to V$ such that $\tilde f(v_j)=v_{f(j)}$, $\forall j\in \{0,\ldots n\}$.
The functor $\mu: \Lambda\longrightarrow\fin$ of \eqref{under} restricts, on the simplicial subcategory $\Delta\subset \Lambda$, to the natural action of $\Delta$ on finite sets that is implemented in the definition of the simplicial structure of $\underline \Delta$. This shows that the latter structure extends to a cyclic structure such that the action is affine and on the vertices is given by
\begin{equation}\label{cycstr}
   \tilde f(v_j)=v_{\mu(f)(j)}\qqq f\in \Hom_\Lambda([n],[m]).
\end{equation}
Note that by construction the cyclic structure \eqref{cycstr} extends to the category $\fin$.
}\end{example}

\section{The geometric realization of a simplicial set and its cyclic structure}

This section reviews and develops on the concept of geometric realization of a simplicial and a cyclic set.

\subsection{The essential geometric morphism  $\hat\Delta\longrightarrow\hat \Lambda$}
The geometric morphism of topoi $\hat\Delta\longrightarrow\hat \Lambda$ associated to the inclusion functor $j:\Delta\hookrightarrow \Lambda$ as in \eqref{embedj}, is an {\em essential geometric morphism} (\cf\cite{MM}, VII.2 Theorem 2). The pullback functor $j^*$ takes a contravariant functor $F: \Lambda\longrightarrow\Se$ to its restriction $F\circ j: \Delta\longrightarrow \Se$. The functor $j^*$ has a left adjoint  $j_!$ which is defined by inducing from $\Delta$ to $\Lambda$
\begin{equation*}
j_!:\hat\Delta\longrightarrow \hat\Lambda\qquad  j_!(S)=S\otimes_\Delta \Lambda \qqq S\in obj(\Se^{\Delta^{\rm op}}).
\end{equation*}
Here, $S$ denotes a simplicial set and the contravariant action of $\Delta$ is viewed as a right action. Next, we use Proposition \ref{decolam} and the decomposition $\Lambda=\Delta {\bf C}$ we referred to in \S\ref{cyccat}, Proposition~\ref{prescyclc}. For $\gamma \in \Aut_\Lambda([n])$ and $\varphi\in \Hom_\Delta([m],[n])$, the composite map $\gamma\circ \varphi$  has a {\em unique decomposition}
\begin{equation}\label{deccyc}
    \gamma\circ \varphi=\gamma_*(\varphi)\circ \varphi^*(\gamma), \qquad
    \gamma_*(\varphi)\in \Hom_\Delta([m],[n]), \ \ \varphi^*(\gamma)\in \Aut_\Lambda([m])
\end{equation}
that encodes the structure of crossed simplicial group of the cyclic category (\cite{Loday} 6.1, Lemma 6.1.5). Moreover, the maps $\gamma_*: \Hom_\Delta([m],[n])\to \Hom_\Delta([m],[n])$ and $\varphi^*: \Aut_\Lambda([n])\to \Aut_\Lambda([m])$ satisfy the following rules  (\opcit 6.1, Proposition 6.1.6)
\begin{equation*}
    (\gamma\gamma')_*(\varphi)=\gamma_*(\gamma'_*(\varphi)), \qquad
    (\varphi\varphi')^*(\gamma)=\varphi'^*(\varphi^*(\gamma)).
\end{equation*}
The decomposition $\Lambda=\Delta {\bf C}$ also provides the following description of the simplicial structure of $j_!(S)$
\begin{equation}\label{jsx}
    j_!(S)_n=S_n\times\Aut_\Lambda([n]), \qquad (x,\gamma)\varphi=(x\gamma_*(\varphi),\varphi^*(\gamma))
\end{equation}
where we encode in the simplicial structure the right action of $\Delta$. The cyclic structure of
$j_!(S)$ is given by
\begin{equation}\label{actgam}
   (x,\gamma)\gamma'=(x,\gamma\gamma')\qquad\forall\gamma,  \gamma'\in \Aut_\Lambda([n]).
\end{equation}
The following identities follow from the crossed simplicial structure and can be derived directly using the uniqueness of the decomposition $\Lambda=\Delta {\bf C}$
\begin{align*}
   \gamma_*(\varphi\varphi')&=\gamma_*(\varphi)(\varphi^*(\gamma))_*(\varphi')\\~\notag\\
\label{crsg1}
   \varphi^*(\gamma\gamma')&=(\gamma'_*(\varphi))^*(\gamma)\varphi^*(\gamma').
\end{align*}
These equalities are used to check the compatibility of the cyclic action \eqref{actgam} with the simplicial structure \eqref{jsx}.
By applying the description of $j_!(S)$ as in \eqref{jsx} to the trivial simplicial set $S=\{pt\}$ one obtains the following 
\begin{prop}
The following formula defines an isomorphism of the  cyclic set $C=j_!(\{pt\})$ with the cyclic set $\yon([0])$ associated to the object $[0]$ of $\Lambda$ by means of the Yoneda functor $\yon: \Lambda\longrightarrow \hat\Lambda$ of \eqref{yoneda}
\begin{equation*}
i: C\stackrel{\sim}{\longrightarrow} \yon([0]),\qquad   (j_!(\{pt\}))_n=\Aut_\Lambda([n])\ni\gamma\mapsto i(\gamma):= f_n\circ \gamma\in \Hom_\Lambda([n],[0])
\end{equation*}
where $f_n$ is the unique element of $\Hom_\Delta([n],[0])$.
\end{prop}
\proof One sees that for $\varphi\in\Hom_\Delta([m],[n])$ the following equalities hold
\begin{equation*}
    i(\varphi^*(\gamma))=f_m\circ \varphi^*(\gamma)=f_n\circ \gamma_*(\varphi)\circ \varphi^*(\gamma)
    =f_n\circ \gamma\circ \varphi=i(\gamma)\circ \varphi
\end{equation*}
The compatibility with the cyclic structure is easily checked. \endproof

\subsection{The geometric realization of a simplicial set}
The geometric realization of a simplicial set $S$ is defined as the quotient space
\begin{equation*}
    |S|:=\left(\coprod_{n\ge 0} (S_n\times \underline\Delta^n)\right)/\sim
\end{equation*}
for the equivalence relation
\begin{equation}\label{equivrel}
    (x\varphi ,t)\sim (x, \varphi_*t)\qquad{\rm for}~ x\in S_n, ~t\in \underline \Delta^m, ~\varphi\in\Hom_\Delta([m],[n]).
\end{equation}
$|S|$ is endowed with the quotient topology and as a set  is given by the tensor product
\begin{equation*}
    |S|=S\otimes_\Delta\underline \Delta=S\otimes_\Delta p_{[0,1]}=(p_{[0,1]})^*(S)
\end{equation*}
where $p_{[0,1]}$ is the point of the topos $\hat \Delta$ associated to the interval $[0,1]$ and $(p_{[0,1]})^*$ is the pullback part of the corresponding geometric morphism $\Se\longrightarrow \hat\Delta$ as described in \eqref{tens}.\vspace{.05in}

The geometric realization of the underlying simplicial set of the cyclic set $C=j_!(\{pt\})$ is the circle: $|C|=\R/\Z$ obtained by gluing together the two endpoints of the interval $I=[0,1]$. More generally one has the following result
\begin{prop}\label{geomrealc}
Let $I=[b,t]$ be an interval and let $p_I$ be the associated point of the topos $\hat \Delta$. Then the map
\begin{equation}\label{ctensdelta}
 \iota_I:C\otimes_\Delta p_I\to I/\sim\qquad   (\tau^a, w)\mapsto \iota_I(\tau^a, w):= w(a)\in I\qqq a\in \{0,\ldots ,n\}
\end{equation}
is a bijection of sets where $C=j_!(\{pt\})$ and $\sim$ is the equivalence relation which identifies the end points  $b$ and $t$ of $I$.
\end{prop}
\proof One has $C=\yon([0])$. The maps $\alpha\in\Hom_\Lambda([n],[0])$ are parameterized by the cosets $a\in \Z/(n+1)\Z$ and  are described, for $f_n\in \Hom_\Delta([n],[0])$,  by
\begin{equation*}
   \alpha_a=f_n\circ \tau^a, \qquad \alpha_a(x)=E(\frac{x-a}{n+1})\qquad x\in \Z
\end{equation*}
where $E(y)$ denotes the integral part of $y\in \R$. Let $\alpha=f_1\circ \tau$ be the map in
$\Hom_\Lambda([1],[0])$ distinct from $f_1$. Then, for $n>1$ and $j\in \{0,\ldots,n-1\}$ one has
\begin{equation}\label{elts}
   \alpha_{j+1}=\alpha\circ \sigma_0\circ\cdots \circ \sigma_{j-1}\circ\sigma_{j+1}\circ \cdots \circ \sigma_{n-1}.
\end{equation}
Notice that $\sigma_j$ is the only omitted term in \eqref{elts} where the $\sigma_i$'s appear following the  increasing indexing order.
By implementing the equivalence relation \eqref{equivrel}, this shows that in $C\otimes_\Delta p_I$, and for any $w\in \Hom_\geq (n^*,I)$, one has
\begin{equation*}
    (\alpha_{j+1},w)\sim (\alpha, (\sigma_0\circ\cdots \circ \sigma_{j-1}\circ\sigma_{j+1}\circ \cdots \circ \sigma_{n-1})_*w)=(\alpha,w\circ s_{n-1}\circ \cdots s_{j+1}\circ s_{j-1}\circ \cdots \circ s_0)
\end{equation*}
where $(w\circ s_{n-1}\circ \cdots s_{j+1}\circ s_{j-1}\circ \cdots \circ s_0)\in \Hom_\geq (1^*,I)=I$. Moreover one also sees that
\begin{equation*}
    (s_{n-1}\circ \cdots s_{j+1}\circ s_{j-1}\circ \cdots \circ s_0)(1)=j+1\in n^*\qqq j\in \{0,\ldots,n-1\}
\end{equation*}
so that
\begin{equation*}
    (\alpha_{j+1},w)\sim (\alpha, w(j+1))\qquad \forall w\in \Hom_\geq (n^*,I).
\end{equation*}
Thus the elements of $C\otimes_\Delta p_I$ are  equivalent to elements either of the form $(\alpha,u)$ for some $u\in I$, or of the form $(\alpha_0,w)$ for some  $w\in \Hom_\geq (n^*,I)$. In this latter case one has
\begin{equation*}
   (\alpha_0,w)\sim (f_0, t_0) , \qquad  t_0\in \Hom_\geq (0^*,I).
\end{equation*}
Then one concludes that $C\otimes_\Delta p_I$ is obtained by gluing the base point $*=(f_0, t_0)$ at the two endpoints of the interval $I\sim\{(\alpha,u)\mid u\in I\}$ and that the map
$\iota_I$ of \eqref{ctensdelta} gives the required bijection. \endproof

\subsection{The geometric realization of a cyclic set}
Next, we re-state Theorem 7.1.4 in \S7.1 of \cite{Loday} in a form that is more suitable to be extended to arbitrary cyclic structures on points of the topos of simplicial sets: we refer to the following sections \S\ref{lawgen} and \S\ref{sectcycpoints}  for precise statements. The notations are the same as in the last (sub)section.
\begin{thm}\label{thmgeom}
$(i)$~For $C=j_!(\{pt\})$, there exists a unique group law on $|C|$ such that
\begin{equation}\label{gplaw}
    (\gamma,u)\cdot(\gamma',\gamma_* u)=(\gamma \gamma',u)\qqq \gamma,\gamma'\in
    \Aut_\Lambda([n]), \ u \in \underline \Delta^n.
\end{equation}
$(ii)$~The map $\iota_{[0,1]}:|C|\to \R/\Z$ as in \eqref{ctensdelta} is a group isomorphism.\vspace{0.05 in}

$(iii)$~Let $S$ be a cyclic set. There exists a unique right action of the group $|C|$ on $|S|$ such that
\begin{equation*}
   (x,u)\triangleleft  (\gamma,u)=(x\gamma^{-1},\gamma_* u) \qqq \gamma\in
    \Aut_\Lambda([n]), \ x\in S_n,\  u \in \underline \Delta^n.
\end{equation*}
\end{thm}
\proof We endow $|C|$ with the group law induced on $|C|$ by the bijection $\iota_{[0,1]}$ of $(ii)$ and show that it fulfills \eqref{gplaw}. Let $u \in \underline \Delta^n$ be an element of the standard simplex, given (the notations as in  \eqref{simplex}) by
\begin{equation*}
    u=\sum_{j=0}^n u_jv_j\in \underline \Delta^n\, , \ \ u=(u_0,\ldots, u_n),\ u_j\geq 0, \ \sum_{j=0}^n u_j=1.
\end{equation*}
Then the map $\sigma^{-1}(u)\in \Hom_\geq(n^*,[0,1])$ associated  to $u$ by \eqref{simplex} is given by
\begin{equation*}
    \sigma^{-1}(u)=\beta, \ \  \beta_a=\sum_{\ell=0}^{a-1} u_\ell,\qquad \forall a\in n^*.
\end{equation*}
Let $\gamma \in \Aut_\Lambda([n])$ and $u \in \underline \Delta^n$: one has $\gamma=\tau^a$ for some $a\in\{0,\ldots ,n\}$, and by \eqref{ctensdelta}
\begin{equation*}
  \iota_{[0,1]}(\tau^a,u)=\beta_a= \sum_{\ell=0}^{a-1} u_\ell\in \R/\Z.
\end{equation*}
Thus one gets
\begin{equation*}
   \iota_{[0,1]}(\tau^b,\tau^a_*u)\sim \sum_{\ell=0}^{b-1} u_{a+\ell}=\sum_{\ell=a}^{a+b-1} u_\ell.
\end{equation*}
Using the relation $\sum_{j=0}^n u_j=1$ one concludes  that the two sides of the following formula  agree modulo $1$:
\begin{equation*}
    \iota_{[0,1]}(\gamma,u)+\iota_{[0,1]}(\gamma',\gamma_* u)=\iota_{[0,1]}(\gamma \gamma',u)\qqq \gamma,\gamma'\in
    \Aut_\Lambda([n]), \ u \in \underline \Delta^n
\end{equation*}
and this proves $(i)$ and $(ii)$.

$(iii)$~First notice that the equality $|C\times S|=|C|\times |S|$ derives from the left exactness of the functor $p_I^*$ for any interval $I$. Moreover we claim that the map
\begin{equation*}
g:  |C\times S|\longrightarrow |S|,\qquad   g(\gamma,x,u)=(x\gamma^{-1},\gamma_* u) \qqq \gamma\in
    \Aut_\Lambda([n]), \ x\in S_n,\  u \in \underline \Delta^n
\end{equation*}
is well defined. This can be checked directly by verifying that, for $\varphi\in \Hom_\Delta([n],[m])$
\begin{equation*}
    g((\gamma,x)\varphi,u)=g(\gamma,x,\varphi_*u)\qqq \gamma\in \Aut_\Lambda([m]), \ x\in S_m, \ u\in
    \underline \Delta^n.
\end{equation*}
One has $(\gamma,x)\varphi=(\varphi^*(\gamma),x\varphi)$ and using \eqref{deccyc} one obtains
\begin{equation*}
    g((\gamma,x)\varphi,u)=\left((x\varphi)\varphi^*(\gamma)^{-1},\varphi^*(\gamma)_*u\right)=
    \left(x\gamma^{-1}\gamma_*(\varphi),\varphi^*(\gamma)_*u\right)
\end{equation*}
since the equality $\gamma\varphi=\gamma_*(\varphi)\varphi^*(\gamma)$ implies that
$\varphi\varphi^*(\gamma)^{-1}=\gamma^{-1}\gamma_*(\varphi)$. Thus one gets
\begin{equation*}
   g((\gamma,x)\varphi,u)=\left(x\gamma^{-1},\gamma_*(\varphi)\varphi^*(\gamma)_*u\right)=
   \left(x\gamma^{-1},(\gamma\varphi)_*u\right)=g(\gamma,x,\varphi_*u).
\end{equation*}
It remains to check that the above map $g: |C|\times |S| \to |S|$ is a right group action. Let $\gamma_j\in \Aut_\Lambda([n])$,  $x\in S_n$, and $u\in
    \underline \Delta^n$. One has by construction
\begin{equation*}
 (x,u)\triangleleft   (\gamma_1,u)=g(\gamma_1,x,u)=(x\gamma_1^{-1},(\gamma_1)_*u)
\end{equation*}
\begin{equation*}
    \left((x,u)\triangleleft   (\gamma_1,u)\right)\triangleleft  (\gamma_2,(\gamma_1)_*u) =
    g(\gamma_2,x\gamma_1^{-1}, (\gamma_1)_*u)=(x\gamma_1^{-1}\gamma_2^{-1},(\gamma_2)_*(\gamma_1)_*u)
    =(x,u)\triangleleft (\gamma_1\gamma_2,u)
\end{equation*}
Together with \eqref{gplaw} this shows that $g$ defines a right group action of the group $|C|$ on $|S|$. \endproof

\section{The cyclic structure associated to an ordered group and a central element}

In this section we show how to associate a cyclic structure (\cf~Definition \ref{defncycpoints}) to a triple $(G,P,z)$, where the pair $(G,P)$ is  a {\em left-ordered group} $G$ (\cf~Definition \ref{defnog} below) and $z\in G $ is  a {\em central positive} element.

\subsection{Ordered groups}

We first recall the definition of a left-ordered group as in \cite{glass} I, \S 1.

\begin{defn}\label{defnog}
A left order on a group $G$ is a semigroup $P\subset G$ (\ie $PP\subseteq P$) such that
$P\cap P^{-1}=\{1\}$ and $P\cup P^{-1}=G$.
\end{defn}
The associated left order on $G$ is the total order $\le$ defined by the relation ($a,b\in G$)
\begin{equation*}
    a\leq b\iff a^{-1}b\in P.
\end{equation*}
By construction the left order $\le$ on $G$ determined by $P$ is invariant by left translations.

\subsection{The cyclic structure associated to $(G,P,z)$}

Next we show how to associate a cyclic structure to a left ordered group $(G,P)$ with a fixed central element $z\in P$. In the following we shall use the notations introduced in \S 2.

\begin{prop}\label{propcp}
Let $(G,P,z)$ be a left ordered group with a fixed central element $z\in P$. The following equality defines a cyclic structure on the point $p_I\in \hat\Delta$ (of the  topos of simplicial sets) associated to the interval $I=[1,z]$
\begin{equation}\label{constructcp}
    \tau(\beta)_j:=\beta_1^{-1}\beta_{j+1} \qqq j\in \{0,\ldots,n\},
    \ \ \tau(\beta)_{n+1}:=z\qqq \beta\in \Hom_\geq(n^*,[1,z]).
\end{equation}
\end{prop}
\proof To an element $\beta\in \Hom_\geq(n^*,[1,z])$ one associates the sequence $\pi(\beta)=\{\pi(\beta)_j\}$ of elements of $P\subset G$ given by 
\begin{equation}\label{pibeta}
    \pi(\beta)_j=\beta_j^{-1}\beta_{j+1} \qqq j\in \{0,\ldots,n\}.
\end{equation}
With $g_j=\pi(\beta)_j$, one has $g_0g_1\cdots g_n=z$. The map $\pi$ determines a bijection of $\cF_{[1,z]}(n^*)=\Hom_\geq(n^*,[1,z])$ with
\begin{equation}\label{fofn}
  F(n):=\{(g_i)_{0\leq i\leq n}\mid g_i\in P, \ g_0g_1\cdots g_n=z\}.
\end{equation}
The inverse map is given by
\begin{equation}\label{betamap}
    \pi^{-1}((g_i))=\beta,\quad \beta_j:= g_0\cdots g_{j-1}\qqq j\in \{0,\ldots,n+1\}.
\end{equation}
Under this change of variables, the covariant action of $\Delta$ is expressed as follows
\begin{align*}
   F( \sigma_j)&:F(n+1)\to F(n),\qquad F(\sigma_j)(g)=(g_0,\ldots, g_{j-1}, g_jg_{j+1},\ldots,g_{n+1})\\
\label{sjdj4}
    F(\delta_j)&:F(n-1)\to F(n), \qquad F(\delta_j)(g)=(g_0,\ldots, g_{j-1},1, g_j,\ldots,g_{n-1}).
\end{align*}
Moreover the action of $\tau$ given by \eqref{constructcp} now reads as follows
\begin{equation}\label{constructcp1}
    F(\tau)(g)_j:=g_{j+1} \qqq j\in \{0,\ldots,n-1\},
    \qquad F(\tau)(g)_n:=g_{0}\qqq g\in F(n).
\end{equation}
Indeed, one has $$\pi(\tau(\beta))_j=\tau(\beta)_j^{-1}\tau(\beta)_{j+1}=\beta_{j+1}^{-1}\beta_{j+2}=\pi(\beta)_{j+1}\qquad \forall j\in \{0,\ldots,n-1\}.$$ For $j=n$: 
$\pi(\tau(\beta))_n=\tau(\beta)_n^{-1}\tau(\beta)_{n+1}=\tau(\beta)_n^{-1}z=z^{-1}\beta_1 z=\beta_1$, while it follows from \eqref{pibeta} that $\pi(\beta)_0=\beta_1$. This shows \eqref{constructcp1}. Note that the cyclic permutation \eqref{constructcp1} preserves the condition $g_0g_1\cdots g_n= z$ in \eqref{fofn} that defines $F(n)$ because the element $z\in G$ is central.

It follows from \eqref{constructcp1} that  $F(\tau)^{n+1}=id$. It remains to check the relations \eqref{pres1} and \eqref{pres2} of  Proposition \ref{prescyclc} which define the presentation of the cyclic category. One has
\begin{equation*}
    F(\tau)\circ F(\sigma_0)(g)=(g_2, \ldots, g_{n+1},g_0g_1)=F(\sigma_n)(g_2, \ldots, g_{n+1},g_0,g_1)=F(\sigma_n)\circ F(\tau)^2 (g).
\end{equation*}
For $j\in \{1, \ldots, n\}$ one has:
\begin{equation*}
   F(\tau)\circ F(\sigma_j)(g)=(g_1, \ldots,g_j g_{j+1},\ldots,g_{n+1},g_0 )=F(\sigma_{j-1})\circ F(\tau)(g).
\end{equation*}
Similarly one has
\begin{equation*}
    F(\tau)\circ F(\delta_0)(g)=F(\tau)(1,g_0, \ldots, g_{n-1})=
    (g_0,\ldots, g_{n-1},1)=
    F(\delta_n) (g)
\end{equation*}
and for $j\in \{1, \ldots, n\}$ one has:
\begin{equation*}
   F(\tau)\circ F(\delta_j)(g)=(g_1, \ldots, g_{j-1},1,g_j,\ldots,g_{n-1},g_0 )=F(\delta_{j-1})\circ F(\tau)(g).
\end{equation*}
This proves that \eqref{constructcp} gives a cyclic structure on the point $p_I$ of $\hat\Delta$. \endproof

The above construction only involves the following subgroup of the group $G$,
\begin{equation}\label{remgprime}
G'=\{g\in G\mid \exists n\in \N, \ z^{-n}\leq g \leq z^n\}
\end{equation}
For any $g\in G'$ there is a largest $k\in \Z$ such that $g\geq z^k$ and one then has $z^{-k}g\in [1,z)$. Thus $G'$ is the subgroup of $G$ generated by the interval $I=[1,z]$ and any of its elements is uniquely of the form $g=z^k v=v z^k$, $k\in \Z$, $v\in [1,z)$. Moreover one has a canonical bijection 
\begin{equation}\label{remgprimebis}
G'/z^\Z=[1,z]/\sim 
\end{equation}
where $\sim$ is the equivalence relation which identifies the end points  $1$ and $z$ of $I=[1,z]$.

\begin{example}\label{eximp}{\rm We give an example (\cf \cite{boy}) of left ordered group $(G,P)$ and central element $z$ such that the group $G'$ is noncommutative. We let $(x_n), n\in \N$ be a dense sequence of real numbers $x_n\in \R$ and consider the left order on the group ${\rm Homeo}_+(\R)$ of order preserving homeomorphisms of $\R$ which is defined by the semigroup
\begin{equation*}
P:=\{ \phi \in {\rm Homeo}_+(\R)\mid \exists n\in \N, \phi(x_n)>x_n, \  \phi(x_j)=x_j  \ \forall j<n\} 
\end{equation*}
Let then $G$ be the subgroup of ${\rm Homeo}_+(\R)$ which is the centralizer of the translation by $1$: $z(x):=x+1, \forall x\in \R$. Thus 
$G=\{ \phi \in {\rm Homeo}_+(\R)\mid  \phi(x+1)=\phi(x)+1 \ \forall x\in \R\}$. One checks that with the induced order, the group $G$ is left ordered and is noncommutative, that $z$ is a central element and that $G'=G$  with $G'$ defined by \eqref{remgprime}.

}\end{example}

\subsection{The action of $G'/z^\Z$ on $p^*S$ for a cyclic set $S$}\label{lawgen}

Let $(G,P,z)$ be a left-ordered group with a fixed central element $z\in G$. We denote by $F$ the extension to the cyclic category (defined in Proposition \ref{propcp})  of the functor $\cF_{[1,z]}: \Delta^{{\rm op}} \longrightarrow \Se$. For a simplicial set $S$ the pullback functor $p_{[1,z]}^*$ applied to $S$ defines by construction the set
\begin{equation*}
    |S|_p=S\otimes_\Delta p_{[1,z]}=(p_{[1,z]})^*(S).
\end{equation*}

We let $G'$ be as in \eqref{remgprime}.
 We can now state the generalization of Theorem \ref{thmgeom} in the presence of  the cyclic structure of Proposition \ref{propcp} on the point $p_{[1,z]}$.

\begin{thm}\label{thmgeom1}
$(i)$~There exists a unique group law on $|C|_p$,  ($C = j_!(\{{pt}\})$), which is given by
\begin{equation*}
    (\gamma,g)\cdot(\gamma',\gamma_* g)=(\gamma \gamma',g)\qqq \gamma,\gamma'\in
    \Aut_\Lambda([n]), \ g \in F(n).
\end{equation*}
$(ii)$~The map $\iota_{[1,z]}: |C|_p\stackrel{\sim}{\to}G'/z^\Z$ is a group isomorphism.\vspace{0.05 in}

$(iii)$~Let $S$ be a cyclic set. There exists a unique right action of the group $|C|_p$ on $|S|_p$ given by
\begin{equation*}
    (x,g)\triangleleft (\gamma,g)=(x\gamma^{-1},\gamma_* g) \qqq \gamma\in
    \Aut_\Lambda([n]), \ x\in S_n,\  g \in  F(n).
\end{equation*}
\end{thm}
\proof By Proposition  \ref{geomrealc} applied to the interval 
$I=[1,z]$  the map
\begin{equation*}
 \iota_{[1,z]}:|C|_p=C\otimes_\Delta p_{[1,z]}\to [1,z]/\sim\qquad   (\tau^a, w)\mapsto \iota_{[1,z]}(\tau^a, w):= w(a)\in I\qqq a\in \{0,\ldots ,n\}
\end{equation*}
is a bijection of sets where $C=j_!(\{pt\})$ and $\sim$ is the equivalence relation which identifies the end points  $1$ and $z$ of $I={[1,z]}$.
As in the proof of Theorem \ref{thmgeom}, we use $ \iota_{[1,z]}$ and  \eqref{remgprimebis} to identify $|C|_p$ with the quotient group $G'/z^\Z$. Thus it is enough, in order to prove $(i)$ and $(ii)$
to show that
\begin{equation*}
  \iota_{[1,z]}(\gamma,g).\iota_{[1,z]}(\gamma',\gamma_* g)=\iota_{[1,z]}(\gamma \gamma',g)\qqq \gamma,\gamma'\in
    \Aut_\Lambda([n]), \ g \in F(n).
\end{equation*}
This follows from \eqref{ctensdelta} and \eqref{betamap} which imply
\begin{equation*}
  \iota_{[1,z]}(\tau^a,g)=\beta_a= g_0\cdots g_{a-1},\qquad  \iota_{[1,z]}(\tau^b,\tau^a_*g)=  g_{a}\cdots g_{a+b-1}.
\end{equation*}
The proof of $(iii)$ is the same as in Theorem \ref{thmgeom}. \endproof

\section{Classification of the cyclic structures}\label{sectcycpoints}

In this section we state and prove the main result of this paper that provides a precise description of cyclic structures. In \S \ref{the last} we relate the notion of cyclic structure with the notion of point of the topos of cyclic sets as in \cite{Moerdjik}.

\subsection{The variety of cyclic structures}

In the following we show that given a point $p=p_I$ of the topos $\hat \Delta$ of simplicial sets, the cyclic structures on $p$ are classified by the (not necessarily commutative) {\em left-ordered group structures} on the ordered set
\begin{equation}\label{ordset}
    G=\left(\Z\times I\right)/\sim \
\end{equation}
Here $\Z\times I$ is endowed with the following lexicographic order of the product  ($I$ has smallest element $b$ and largest element $t$):
\begin{equation}\label{lex}
   n>m\implies (n,u)> (m,v), \qquad (n,u)\geq (n,v)\iff u\geq v\qquad\quad u,v\in I.
\end{equation}
The equivalence relation in \eqref{ordset} identifies $(n,t)\sim (n+1,b)$ for all $n\in \Z$. By construction there is a natural injection of sets
\begin{equation}\label{inject}
   c: \Z\hookrightarrow G, \qquad c(n)=(n,b)\qquad n\in \Z.
\end{equation}
The main result of this paper is the following
\begin{thm}\label{main} Let $I$ be an interval and let $p_I$ be the associated point of the topos $\hat \Delta$ of simplicial sets. Then a cyclic structure on $p_I$ corresponds to a group law on $ G=\left(\Z\times I\right)/\sim$~ such that:
\begin{enumerate}
  \item  The order relation on $G$ is left invariant
   \item The restriction of the group law of $G$ on $c(\Z)\times G$  is commutative and is given by
   \[
   c(n)(m,u)=(m,u) c(n)=(n+m,u),\qquad \forall n,m\in \Z, ~u\in I.
   \]
\end{enumerate}
\end{thm}
Note that (ii) implies that the  map $c:\Z\hookrightarrow G$ of \eqref{inject} is a group homomorphism whose range is contained in the center of $G$.
Before to start the proof of Theorem \ref{main} we state some simple corollaries and provide some examples showing the variety of cyclic structures one can define on a point $p=p_I$ of $\hat\Delta$.

\begin{cor}\label{hom}

$(i)$~Let $I$ be an interval and let $p_I$ be the associated point of $\hat \Delta$. If $p_I$ admits a cyclic structure then the ordered set \eqref{ordset} is homogeneous.\vspace{0.01in}

$(ii)$~There exist points of $\hat \Delta$ that cannot support any cyclic structure.\vspace{0.01in}

$(iii)$~The points of $\hat \Delta$ associated to {\em finite} intervals have a unique cyclic structure.
\end{cor}
\proof $(i)$~follows from $(i)$ of Theorem \ref{main} since left translations by $G$ respect the order and act transitively.

$(ii)$~Let $I$ be an interval with a unique non-isolated point $u\in (b,t)$ for the order topology (\ie the topology with a basis given by open subintervals);  all the other points are isolated \ie open. Then the ordered set \eqref{ordset} is not homogeneous thus $(i)$ implies that the point $p_I$  cannot have any cyclic structure.

$(iii)$~A point $p_I$ with $I$ a finite interval corresponds to the interval $n^*$, for some $n\geq 0$. The associated ordered set \eqref{ordset} is isomorphic to $\Z$ with the usual order, and the inclusion $c:\Z\hookrightarrow \Z$ of \eqref{inject} is defined by multiplication by $n+1$. Then, the existence of a group law fulfilling the conditions of Theorem \ref{main} is immediate. Conversely, let $e=(0,0)$ be the neutral element of the group $G$ and let $f=(0,1)$ be the subsequent element after $e$ for the order relation \eqref{lex}. We write the group law of $G$ multiplicatively. The sequence $\{f^n\}$ is increasing: in fact $f>e$, thus by applying the left invariance property $(i)$ of Theorem \ref{main} one gets  $f^{n+1}>f^n$.
Moreover the open intervals $(f^n,f^{n+1})$ are empty since
\begin{equation*}
    f^n<x<f^{n+1}\implies e<f^{-n}x<f.
\end{equation*}
It follows that $f^{n+1}=(1,0)=c(1)$ with $f$ generating the group $G$ with positive part provided by the $f^n$'s for $n\geq 0$. This proves the required uniqueness.
\endproof

\begin{cor}\label{hom1}
$(i)$~Let $I$ be an interval and let $\Gamma=\Aut(I)$ be the group  of  order automorphisms  of $I$. Then $\Gamma$ acts on the set of  cyclic structures on  the associated point $p_I$ of $\hat \Delta$.\vspace{0.01in}

$(ii)$~Let $I=[0,1]\subset \R$, then the associated point $p_{[0,1]}$ of  $\hat \Delta$ admits a continuum of cyclic structures.
\end{cor}
\proof $(i)$~The two properties fulfilled by a cyclic structure on $p_I$ as in  Theorem \ref{main} remain true if one conjugates a group law $(g,g')\mapsto gg'$ by an order automorphism $\psi$, \ie if one  replaces the group law by the following modified one
\begin{equation*}
    (g,g')\mapsto \psi(\psi^{-1}(g)\psi^{-1}(g')).
\end{equation*}
$(ii)$~The group $\Aut([0,1])$ of order automorphisms of the interval $[0,1]$ acts freely on the set of cyclic structures on $p_{[0,1]}$ which are obtained from the usual ordered group structure on $\R$.\endproof

\subsection{Proof of Theorem \ref{main}}\label{mainsect}

Let $p$ be a point of $\hat \Delta$  and let $\cF=p^*\circ \yon: \Delta \longrightarrow \Se$ be the associated (filtering) functor. Let $\tilde \cF:\Lambda \longrightarrow \Se$ be an extension of $\cF$ to the cyclic category that defines a cyclic structure on $p$. Let $I$ be the interval associated to $\cF$ so that $\cF=\cF_I$ (\cf~\eqref{fsubi}). The action of $\Aut_\Lambda([n])$ arising from the extension $\tilde\cF$ is denoted as follows
\begin{equation*}
    \tilde \cF(\gamma)(u)= \gamma_* u \in \cF(n)=\Hom_\leq(n^*,I)\qqq \gamma \in
    \Aut_\Lambda([n]), \ u \in \Hom_\leq(n^*,I).
\end{equation*}
More generally we shall use the notation
\begin{equation*}
\tilde \cF(\phi)(u)= \phi_* u \in \cF(m)=\Hom_\leq(m^*,I)\qqq \phi \in
    \Hom_\Lambda([n],[m]), \ u \in \Hom_\leq(n^*,I).
\end{equation*}
The pullback functor $p_I^*$ associates to any simplicial set $S$ the set (\cf\eqref{gerea})
\begin{equation*}
   (p_{I})^*(S)= |S|_p=S\otimes_\Delta p_I.
\end{equation*}
By Proposition \ref{geomrealc} and for $C = j_!(\{pt\})=\yon([0])$, one has a bijection $\iota_I:|C|_p\to I/\sim$~, where $\sim$ is the equivalence relation which identifies the end points  $b$ and $t$ of $I$.
Next statement defines  a group law on $|C|_p$
\begin{lem}\label{lempf0} There exists a unique group law on $|C|_p$ such that
\begin{equation}\label{gplaw2}
    (\gamma,u)\cdot(\gamma',\gamma_* u)=(\gamma \gamma',u)\qqq \gamma,\gamma'\in
    \Aut_\Lambda([n]), \ u \in \Hom_\leq(n^*,I).
\end{equation}
\end{lem}
\proof Let $S$ be any cyclic set. There is an equality of geometric realizations $|C\times S|_p=|C|_p\times |S|_p$ that derives, from the left exactness of the functor $p_I^*$.  Moreover  we claim that the map 
\begin{equation*}
g_p:  |C\times S|_p\longrightarrow |S|_p,\quad   g_p(\gamma,x,u)=(x\gamma^{-1},\gamma_* u),~ \forall\gamma\in
    \Aut_\Lambda([n]), \ x\in S_n,\  u \in \Hom_\leq(n^*,I)
\end{equation*}
is well defined. This can be checked directly by verifying that, for $\varphi\in \Hom_\Delta([n],[m])$
\begin{equation*}
    g_p((\gamma,x)\varphi,u)=g_p(\gamma,x,\varphi_*u),~\forall\gamma\in \Aut_\Lambda([m]), \ x\in S_m, \ u\in
    \Hom_\leq(n^*,I).
\end{equation*}
One has $(\gamma,x)\varphi=(\varphi^*(\gamma),x\varphi)$, and using \eqref{deccyc} one obtains
\begin{equation*}
    g_p((\gamma,x)\varphi,u)=\left((x\varphi)\varphi^*(\gamma)^{-1},\varphi^*(\gamma)_*u\right)=
    \left(x\gamma^{-1}\gamma_*(\varphi),\varphi^*(\gamma)_*u\right)
\end{equation*}
since the equality $\gamma\varphi=\gamma_*(\varphi)\varphi^*(\gamma)$ implies 
$\varphi\varphi^*(\gamma)^{-1}=\gamma^{-1}\gamma_*(\varphi)$. Thus one gets
\begin{equation*}
   g_p((\gamma,x)\varphi,u)=\left(x\gamma^{-1},\gamma_*(\varphi)\varphi^*(\gamma)_*u\right)=
   \left(x\gamma^{-1},(\gamma\varphi)_*u\right)=g_p(\gamma,x,\varphi_*u).
\end{equation*}
Applying this result to $S=C$ one derives the map
\begin{equation}\label{grprod}
h_p: |C|_p\times |C|_p\to |C|_p,\quad    h_p\left((\gamma,u),(\gamma',u)  \right) =(\gamma'\gamma^{-1},\gamma_*u),~\forall\gamma,\gamma'\in
    \Aut_\Lambda([n]), \ u \in \Hom_\leq(n^*,I).
\end{equation}
This implies
\begin{equation*}
    h_p\left((\gamma,\gamma^{-1}_*v),(\gamma',\gamma^{-1}_*v)  \right) =(\gamma'\gamma^{-1},v),~\forall \gamma,\gamma'\in
    \Aut_\Lambda([n]), \ v\in \Hom_\leq(n^*,I).
\end{equation*}
Moreover by following the above construction one also shows that there exists a map 
\begin{equation*}
 j:|C|_p\to |C|_p,\qquad   j(\gamma,v)=(\gamma^{-1},\gamma_*v)\qqq \gamma\in
    \Aut_\Lambda([n]), \ v\in \Hom_\leq(n^*,I).
\end{equation*}
Thus the map $k: |C|_p\times |C|_p\to |C|_p$, ~$k(x,y)=h_p(j(x),y)$ fulfills the rule
\begin{equation*}
    k\left((\gamma,v),(\gamma',\gamma_*v)  \right) =(\gamma'\gamma,v)\qqq \gamma,\gamma'\in
    \Aut_\Lambda([n]), \ v\in \Hom_\leq(n^*,I).
\end{equation*}
We  set $x.y=k(x,y)$
and show that this defines a group law on $|C|_p$.
For $j\in\{0,1,2\}$, let $\gamma_j\in \Aut_\Lambda([n])$,   and $u\in
   \Hom_\leq(n^*,I)$. We use the above map $g_p$ with $S=C$ and set
\begin{equation*}
 (\gamma_0,u)\triangleleft   (\gamma_1,u):=h_p\left((\gamma_1,u),(\gamma_0,u)\right)=(\gamma_0\gamma_1^{-1},(\gamma_1)_*u).
\end{equation*}
One also has
\begin{equation*}
    \left((\gamma_0,u)\triangleleft   (\gamma_1,u)\right)\triangleleft  (\gamma_2,(\gamma_1)_*u) =
   (\gamma_0\gamma_1^{-1},(\gamma_1)_*u) \triangleleft  (\gamma_2,(\gamma_1)_*u)\end{equation*}
   \begin{equation*}
   =(\gamma_0\gamma_1^{-1}\gamma_2^{-1},(\gamma_2)_*(\gamma_1)_*u)
    =(\gamma_0,u)\triangleleft (\gamma_1\gamma_2,u).
   \end{equation*}
With the above notations one derives
\begin{equation*}
    (x \triangleleft y)\triangleleft z= x\triangleleft (y.z)\qqq x,y,z\in |C|_p.
\end{equation*}
The elements $(1,u)$ all represent the same element $1\in |C|_p$ and  one has
\begin{equation*}
    1\triangleleft x=j(x)\qqq x\in  |C|_p.
\end{equation*}
The associativity of the product $x.y$ then follows from the following equalities
\begin{equation*}
   1\triangleleft((x.y).z)=1\triangleleft x\triangleleft y \triangleleft z 
   =1\triangleleft(x.(y.z)).
\end{equation*}
Finally, any $x\in |C|_p$ has an inverse $j(x)$, using the fact that $x \triangleleft x=1$.\endproof

In the next step we construct a 2-cocycle $c$ on the group $|C|_p$ with values in the additive group $\Z$. First we introduce the natural triangulation of $|C|_p\times |C|_p$.

\begin{lem}\label{lempf1}
$(i)$~The simplicial structure of $C\times C$ is given by degenerate simplices except for the two $2$-dimensional simplices $T_1=(\tau_2,\tau_2^2)$ and $T_2=(\tau_2^2,\tau_2)$, the three $1$-dimensional simplices $L_1=(1,\tau_1)$, $L_2=(\tau_1,1)$ and $L_3=(\tau_1,\tau_1)$ and the zero dimensional simplex $(1,1)$.

$(ii)$~The face maps $\delta_j$ ($j=0,1,2$) are given by the following table
\begin{equation*}
\begin{array}{ccc}
  T_1\delta_0=L_1, && T_2\delta_0=L_2 \\
  T_1\delta_1=L_3, && T_2\delta_1=L_3 \\
  T_1\delta_2=L_2, && T_2\delta_2=L_1.
\end{array}
\end{equation*}

\end{lem}
\proof $(i)$~There are a priori nine $2$-dimensional simplices in $C\times C$ associated to the pairs
$(\tau_2^a,\tau_2^b)$, for $a,b\in\{0,1,2\}$. To show that with the exception of $T_1=(\tau_2,\tau_2^2)$ and $T_2 = (\tau_2^2,\tau_2)$ the other seven  are degenerate we compute the degeneracy maps $\sigma_i$ ($i=0,1$) from the lower dimensional simplices
\begin{equation*}
    \begin{array}{cccc}
       (1,1)\sigma_0=(1,1), & (1,\tau_1)\sigma_0=(1,\tau_2^2), & (\tau_1,1)\sigma_0=(\tau_2^2,1), & (\tau_1,\tau_1)\sigma_0=(\tau_2^2,\tau_2^2) \\
      (1,1)\sigma_1=(1,1), & (1,\tau_1)\sigma_1=(1,\tau_2), & (\tau_1,1)\sigma_1=(\tau_2,1), &
      (\tau_1,\tau_1)\sigma_1=(\tau_2,\tau_2).
     \end{array}
\end{equation*}
$(ii)$~One obtains the table of degeneracies by direct computation. \endproof

\begin{figure}
\begin{center}
\includegraphics[scale=0.5]{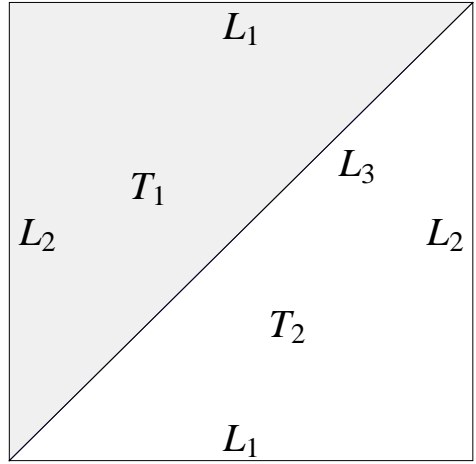}
\caption{The natural triangulation on $|C|_p\times |C|_p$}\label{cyctop}
\end{center}
\end{figure}

Figure \ref{cyctop} shows this natural triangulation on $|C|_p\times |C|_p$.

\begin{rem}
{\rm One can compare the natural triangulation of   $|C|_p\times |C|_p$ with the triangulation of $j_!(C)$ used in \cite{Loday} Lemma 7.1.10, using the following fact.
Let $S$ be a cyclic set, then the following map gives an isomorphism of underlying simplicial sets
\begin{equation*}
   \vartheta_S:j_!(S)\to C\times S, \ \ \vartheta_S(x,\gamma):=(\gamma,x\gamma)\qqq
   \gamma \in\Aut_\Lambda([n]), \ x\in S_n.
\end{equation*}
To prove this statement it is enough to check that  $\vartheta_S$ is compatible with the simplicial structure. One has
\begin{equation*}
   \vartheta_S(x,\gamma)\phi=(\phi^*(\gamma),x\gamma\phi)=
   (\phi^*(\gamma),x\gamma_*(\phi)\phi^*(\gamma))=
   \vartheta_S(x\gamma_*(\phi),\phi^*(\gamma))=\vartheta_S((x,\gamma)\phi).
\end{equation*}
}\end{rem}

Next result uses the above natural  triangulation on $|C|_p\times |C|_p$ to provide an explicit description of the group law on $|C|_p$ and to construct a 2-cocycle $c$ on the group $|C|_p$.

\begin{lem}\label{lempf2}
$(i)$~The group law on $|C|_p=I/\sim$~ is given as follows
\begin{equation}\label{2cocycle0}
    x.y=\left\{
             \begin{array}{ll}
              \left( \tau_2(1,\tau_1(x),y,z)\right)(1) & \hbox{if}\quad\tau_1(x)\leq y \\
               \left( \tau_2^2(1,y,\tau_1(x),z)\right)(2) & \hbox{if}\quad\tau_1(x)> y
             \end{array}
           \qquad \forall x,y\in I\right.
\end{equation}
where $\tau_n$ is the cyclic action on $\cF_I(n)=\Hom_\geq(n^*,I)$.

$(ii)$~The following equality defines a $2$-cocycle $c\in Z^2(|C|_p,\Z)$
normalized by the condition $c(1,x)=c(x,1)=0$, $\forall x\in |C|_p$, where $1\in |C|_p$ is the class of $b\sim t$ and the neutral element of the group law \eqref{2cocycle0}
\begin{equation}\label{2cocycle}
    c(x,y)=\left\{
             \begin{array}{ll}
               0 & \hbox{if}\ \tau_1(x)> y \\
               1 & \hbox{if}\ \tau_1(x)\leq y
             \end{array}
            \qquad \forall x,y\in I\setminus \{b,t\}.\right.
\end{equation}

\end{lem}
\proof $(i)$~On the simplex $T_1=(\tau_2,\tau_2^2)$ the isomorphism $\iota_I\times \iota_I$ of  \eqref{ctensdelta} is given by
\begin{equation*}
  (\tau_2,\tau_2^2,u)=\left((\tau_2,u),(\tau_2^2,u)\right)\stackrel{\iota_I\times \iota_I}{\mapsto} ( u(1),u(2))\qqq u\in \Hom_\geq(2^*,I).
\end{equation*}
Similarly, on the simplex $T_2=(\tau_2^2,\tau_2)$ one has
\begin{equation*}
  (\tau_2^2,\tau_2,u)=\left((\tau_2^2,u),(\tau_2,u)\right)\stackrel{\iota_I\times \iota_I}{\mapsto}  ( u(2),u(1))\qqq u\in \Hom_\geq(2^*,I).
\end{equation*}
By applying in both cases the map  $h_p: |C|_p\times |C|_p\to |C|_p$ as in \eqref{grprod} we get on $T_1$
\begin{equation*}
   (\tau_2,\tau_2^2,u)=\left((\tau_2,u),(\tau_2^2,u)\right)\stackrel{ h_p}{\mapsto} (\tau_2,\tau_2(u))\stackrel{\iota_I}{\mapsto} \tau_2(u)(1)\qqq u\in \Hom_\geq(2^*,I)
\end{equation*}
and on $T_2$
\begin{equation*}
  (\tau_2^2,\tau_2,u)=\left((\tau_2^2,u),(\tau_2,u)\right)\stackrel{ h_p}{\mapsto} (\tau_2^2,\tau_2^2(u))\stackrel{\iota_I}{\mapsto} \tau_2^2(u)(2)\qqq u\in \Hom_\geq(2^*,I).
\end{equation*}
Thus using the definition $k(x,y)=h_p(j(x),y)$ as in the proof of Lemma \ref{lempf0} one obtains \eqref{2cocycle0}.

$(ii)$~We first define a function  $\rho(x,y)$ on $|C|_p\times |C|_p$ with values in $\{0,1\}$ in terms of the triangulation of Lemma \ref{lempf1},
 \begin{equation*}
     \rho=1\ \  \text{on}\ \ T_1\setminus(L_1\cup L_2), \ \ \ \rho=0\ \ \text{elsewhere.}
 \end{equation*}
 Thus $\rho(x,y)=0$ if $x$ or $y$ is the base point $1\in |C|_p$, and for $x,y\in I\setminus\{b,t\}$ one has
 \begin{equation}\label{triang1}
   \rho(x,y)=\left\{
             \begin{array}{ll}
               0, & \hbox{if}\ x> y \\
               1, & \hbox{if}\ x\leq y
             \end{array}\right.
 \end{equation}
 Thus comparing with \eqref{2cocycle} we get $c(x,y)=\rho(j(x),y)$ for all $x,y\in |C|_p$.

 The cocycle equation for a function of two variables $f(a,b)$ on a group $G$ is equivalent to
 \begin{equation*}
    \sum_{j=0}^3(-1)^j\tilde f(w_0,\ldots, \hat w_j, \ldots,w_3)=0 \qqq w_i\in G,
\end{equation*}
where $\tilde f$ is the left invariant function of three variables
\begin{equation*}
    \tilde f(w_0,w_1,w_2):=f(w_0^{-1}w_1,w_1^{-1}w_2) \qqq w_i\in G.
\end{equation*}
For the function $c$ of \eqref{2cocycle} on the group $|C|_p$ we use \eqref{gplaw2}  in the form $(\gamma,u)^{-1}.(\gamma',u)=(\gamma^{-1}\gamma',\gamma_*u)$, and get the equality
\begin{equation*}
    \tilde c\left((\gamma_0,u),(\gamma_1,u),(\gamma_2,u)\right)=
    c\left((\gamma_0^{-1}\gamma_1,(\gamma_0)_*u),(\gamma_1^{-1}\gamma_2,(\gamma_1)_*u)\right)
\end{equation*}
Then we use the equality $c(x,y)=\rho(j(x),y)$ and obtain
\begin{equation}\label{2cc}
 \tilde c\left((\gamma_0,u),(\gamma_1,u),(\gamma_2,u)\right)=
 \rho\left((\gamma_0\gamma_1^{-1},(\gamma_1)_*u),(\gamma_1^{-1}\gamma_2,(\gamma_1)_*u)\right)
\end{equation}
It remains to show that this function of three variables $\tilde c(w_0,w_1,w_2)$ fulfills
\begin{equation}\label{coc3}
    \sum_{j=0}^3(-1)^j\tilde c(w_0,\ldots, \hat w_j, \ldots,w_3)=0.
\end{equation}
 Let $\omega$ be the function of three variables on $|C|_p\simeq I\setminus \{t\}$ given by
\begin{equation*}
    \omega(x,y,z):=\left\{
             \begin{array}{ll}
               0 & \hbox{if}\ x= y \ \hbox{or}\ y=z \\
                0 & \hbox{if}\ x<z \ \hbox{and}\ y\in [x,z] \\
                 0 & \hbox{if}\ z<x \ \hbox{and}\ y\notin [z,x] \\
               1 & \hbox{otherwise}.
             \end{array}\right.
\end{equation*}
By construction  $\omega$  depends only upon the order structure on $I$ and not on the cyclic structure. To verify that it fulfills \eqref{coc3} it is enough to check it in the special case of  $I=[0,1]$ since this interval is universal for finite order configurations in arbitrary intervals. In this case one finds that $\omega$ is the coboundary of the following function of two variables
\begin{equation*}
    \ell(x,y):=\left\{
             \begin{array}{ll}
               y-x & \hbox{if}\ x\leq  y \\
               1-(x-y) & \hbox{if}\ x>y
             \end{array}\right.
\end{equation*}
in other words  in this special case one has 
\begin{equation*}
    \omega(x,y,z)=\ell(y,z)-\ell(x,z)+\ell(x,y).
\end{equation*}
Next we turn back to the general case and show that  $\omega(w_0,w_1,w_2)=\tilde c(w_0,w_1,w_2)$. More precisely we prove the equality on each of the simplices of the canonical triangulation of $|C|_p\times |C|_p\times |C|_p=|C\times C\times C|_p$. We start with the simplices of lowest dimension and parametrize their elements in the form  $(\gamma_0,\gamma_1,\gamma_2,u)$, where for simplices of dimension $n\in \N$ one has $\gamma_j \in \Aut_\Lambda([n])$ and $u\in \Hom_\geq(n^*,I)$. Since we proceed by considering first  the lower dimensional simplices, we can assume that $u$ is non-degenerate. Then it follows immediately that $\gamma_*u$ is also non-degenerate $\forall\gamma\in \Aut_\Lambda([n])$ since an equality such as
$\gamma_*u=\phi_*v
$, for $v\in \Hom_\geq((n-1)^*,I)$ and $\phi\in \Hom_\Delta([n-1],[n])$, implies $u=\psi_*(w)$ for some 
$w\in \Hom_\geq((n-1)^*,I)$ and $\psi\in \Hom_\Delta([n-1],[n])$. To compute $\tilde c(w_0,w_1,w_2)$ for $w_j=(\gamma_j,u)$, we use  \eqref{2cc} and we can then assume that $(\gamma_1)_*u$ is  non-degenerate. Both $\omega(w_0,w_1,w_2)$ and $\tilde c(w_0,w_1,w_2)$ take a constant value on the non-degenerate elements of each simplex and we can then compare these values and check that they are the same. Both $\omega$ and $\tilde c$ take the value $0$ on the $0$-dimensional simplex. We label the simplices of dimension $n$ by the three indices $a_j\in \{0, \ldots,n\}$ corresponding to $\gamma_j=\tau^{a_j}$ and we only consider the non-degenerate simplices.
Next tables list, for higher dimensional simplices, the  values of the three indices $a_j\in \{0, \ldots,n\}$, the value of  $\omega(w_0,w_1,w_2)$ and of the two indices $b_j$ such that $\gamma_0\gamma_1^{-1}=\tau^{b_0}$, $\gamma_2\gamma_1^{-1}=\tau^{b_1}$ and finally the corresponding value of  $\rho\left((\gamma_0\gamma_1^{-1},(\gamma_1)_*u),(\gamma_1^{-1}\gamma_2,(\gamma_1)_*u)\right)$. 
One has $7$ non-degenerate simplices of dimension $1$ obtained using permutations as the following simplices
$(w_{\sigma(0)}=w_{\sigma(1)}=w_{\sigma(2)})$,  $(w_{\sigma(0)}=0,  \ w_{\sigma(1)}=w_{\sigma(2)})$,
$(w_{\sigma(0)}=w_{\sigma(1)}=0)$.
The corresponding table is
$$
\begin{array}{ccc}
 (1, & 0, & 0)\to  0\\
 (0, & 1, & 0)\to  1\\
(0, & 0, & 1)\to  0\\
 (1, & 1, & 0)\to  0\\
 (1, & 0, & 1)\to  1\\
  (0, & 1, & 1)\to  0\\
(1, & 1, & 1)\to  0\\
\end{array},
\begin{array}{cc}
 (1, & 0)\to  0\\
(1, & 1)\to  1\\
(0, & 1)\to  0\\
(0, & 1)\to  0\\
(1, & 1)\to  1\\
(1, & 0)\to  0\\
(0, & 0)\to  0\\
\end{array}
$$ 
One has $12$ non-degenerate simplices of dimension $2$ obtained using permutations as the following simplices
$(w_{\sigma(0)}=w_{\sigma(1)}<w_{\sigma(2)})$,   $(w_{\sigma(0)}=0,  \ w_{\sigma(1)}<w_{\sigma(2)})$.
The corresponding table is
$$
\begin{array}{ccc}
 (0, & 1, & 2)\to  0\\
 (0, & 2, & 1)\to  1\\
(1, & 0, & 2)\to  1\\
 (1, & 1, & 2)\to  0\\
 (1, & 2, & 0)\to  0\\
  (1, & 2, & 1)\to  1\\
(1, & 2, & 2)\to  0\\
(2, & 0, & 1)\to  0\\
(2, & 1, & 0)\to  1\\
 (2, & 1, & 1)\to  0\\
(2, & 1, & 2)\to  1\\
 (2, & 2, & 1)\to  0\\
\end{array},
\begin{array}{cc}
 (2, & 1)\to  0\\
(1, & 2)\to  1\\
(1, & 2)\to  1\\
(0, & 1)\to  0\\
(2, & 1)\to  0\\
(2, & 2)\to  1\\
(2, & 0)\to  0\\
(2, & 1)\to  0\\
(1, & 2)\to  1\\
 (1, & 0)\to  0\\
(1, & 1)\to  1\\
(0, & 2)\to  0\\
\end{array}
$$
One has $6$ non-degenerate simplices of dimension $3$, corresponding to permutations $\sigma$ by the ordering $w_{\sigma(0)}<w_{\sigma(1)}<w_{\sigma(2)}$. 
The corresponding table is
$$
\begin{array}{ccc}
 (1, & 2, & 3)\to  0\\
 (1, & 3, & 2)\to  1\\
(2, & 1, & 3)\to  1\\
 (2, & 3, & 1)\to  0\\
 (3, & 1, & 2)\to  0\\
  (3, & 2, & 1)\to  1\\
\end{array},
\begin{array}{cc}
 (3, & 1)\to  0\\
(2, & 3)\to  1\\
(1, & 2)\to  1\\
(3, & 2)\to  0\\
(2, & 1)\to  0\\
(1, & 3)\to  1\\
\end{array}
$$
This suffices to check the equality  $\omega(w_0,w_1,w_2)=\tilde c(w_0,w_1,w_2)$. In fact the above verifications correspond to the special cases of the (unique) cyclic structures on the intervals $n^*$, for $n\in\{0,1,2,3\}$.
\endproof
Using the $2$-cocycle $c$ defined in Lemma \ref{lempf2} one obtains a group extension of the form
\begin{equation*}
    0\to \Z\to G\to |C|_p\to 1.
\end{equation*}
The elements of $G$ are pairs $(n,x)$, with $n\in \Z$ and $x\in |C|_p=I/\sim$. The group law on $G$ is given by 
\begin{equation*}
   (n,x).(m,y)=(n+m+c(x,y),x.y).
\end{equation*}
The proof of Theorem \ref{main} follows immediately from next statement.  
\begin{lem}
$(i)$~The subset $P=\{(n,x)\mid n\geq 0\}\subset G$ is a semigroup
that satisfies the following properties: $P\cap P^{-1}=\{1\}$, $P\cup P^{-1}=G$.\newline
$(ii)$~The map $\Z\times I\ni (n,u) \mapsto (n,x)\in G$, with $x=$ (class of $u)\in |C|_p$, defines an order isomorphism on the (ordered) set $G = (\Z\times I)/\sim$,~ with the left invariant order on $G$ associated to $P\subset G$.
\end{lem}
\proof $(i)$~By construction of the $2$-cocycle $c$ one has $c(x,y)\geq 0$ for all $x,y\in |C|_p$ thus the subset $P$ is a semigroup. For $x\in |C|_p$, $x\neq 1$, $c(x^{-1},x)=1$, thus the inverse of $(n,x)\in P$ is $(-n-1,x^{-1})$. Then if $n\geq 0$, $(-n-1,x^{-1})\notin P$. This shows that $P\cap P^{-1}=\{1\}$. One sees easily that $P\cup P^{-1}=G$, since for any $n\in \Z$ one has either $n\geq 0$ or $-n-1\geq 0$.\newline
$(ii)$~The left invariant order on $G$ associated to $P\subset G$ is defined as follows
\begin{equation*}
    (n,x)\leq (m,y)\iff (n,x)^{-1}(m,y)\in P.
\end{equation*}
Assuming $x\neq 1$, one has $(n,x)^{-1}=(-n-1,x^{-1})$, hence $(n,x)^{-1}(m,y)=(-1+m-n+c(x^{-1},y), x^{-1}y)$, 
so that (for $x\neq 1$) the above order can be described as
\begin{equation*}
    (n,x)\leq (m,y)\iff -1+m-n+c(x^{-1},y)\geq 0.
\end{equation*}
This corresponds to the lexicographic order since $m>n$ implies $-1+m-n+c(x^{-1},y)\geq 0$ and if $m=n$ one derives using \eqref{triang1} 
\begin{equation*}
    c(x^{-1},y)\geq 1\iff x\leq y.
\end{equation*}
Finally, for $x=1$ one has $(n,1)^{-1}=(-n,1)$ and $ (n,1)^{-1}(m,y)=(m-n,y)$ so that 
\begin{equation*}
  (n,1)\leq (m,y)\iff m-n\geq 0
\end{equation*}
which again corresponds to the lexicographic order. \endproof

\subsection{Relation with the points of the topos of cyclic sets}\label{the last}

In \cite{Moerdjik} it was shown that the topos of cyclic sets is the classifying topos of \emph{abstract circles}. By definition an
abstract circle $C$ is given by the following structure
\begin{equation*}
C=(P,S,\partial_0,\partial_1,0,1,*,\cup)
\end{equation*}
where $P$ and $S$ are sets, $\partial_j:S\to P$ are maps as well as $P\ni x\to 0_x\in S$ and $P\ni x\to 1_x\in S$, $*:S\to S$ is an involution, and $\cup$ is a partially defined map from a subset of $S\times S$ to $S$. The model to keep in mind for such a structure is the one provided by a subset $P\subset S^1$ of the oriented unit circle with $S$ the set of all positively oriented segments with end points in $P$.
The axioms fulfilled by abstract circles are the following:

$1.$~Non-triviality. $P\neq\emptyset$  and for any $x,y\in P$ there exists at least one $a\in S$ such that $\partial_0 a=x$, $\partial_1 a=y$. For any $x\in P$ the segments $0_x$ and $1_x$ are distinct.

$2.$~ Equational. $a^{**}=a$, $\partial_0(a^*)=\partial_1(a)$, $\partial_0( 0_x)=x=\partial_1( 0_x)$, $0_x^*=1_x$, if $\partial_0 a=\partial_1 a=x$ then $a=0_x$ or $a=1_x$.

$3.$~ Concatenation.
\begin{enumerate}
  \item $a\cup b$ exists only if $\partial_1 a=\partial_0 b$ and in that case $\partial_1 (a\cup b)=\partial_1 b$ and $\partial_0 (a\cup b)=\partial_0 a$.
  \item $a\cup b=c \Longleftrightarrow c^*\cup a=b^* $.
  \item If $a\cup b$ and $(a\cup b)\cup c$ exist then so do $b\cup c$ and $a\cup (b\cup c)$ and 
$(a\cup b)\cup c=a\cup (b\cup c)$.
  \item $a\cup b=0_x \implies a=0_x$ 
  \item If $\partial_0 a=x$ then $0_x\cup a=a$.
  \item If $\partial_1 a=\partial_0 b$ then at least one of $a \cup b$ and $b^*\cup a^*$ exist.
\end{enumerate}
We introduce a closely related notion
\begin{defn} An archimedean set is a pair $(X,\theta)$ where
$X$ is a non-empty totally ordered set and $\theta\in \Aut X$ is an order automorphism, with $\theta(x)>x$, $\forall x\in X$ and fulfilling the following archimedean property
\begin{equation}\label{arch}
\forall x,y\in X, \ \exists n\in \N: \quad y\leq \theta^n(x).
\end{equation}
\end{defn}

To any  archimedean set $(X,\theta)$ we associate an abstract circle $C=X/\theta$ as follows
\begin{itemize}
\item $P=X/\sim$~ is the orbit space for the action of $\Z$ on $X$ given by powers of $\theta$.
\item $S$ is  the orbit space for the action of $\Z$ on the set of pairs $(x,y)\in X^2$, with $x\leq y\leq \theta(x)$.
\item $\partial_0 (x,y)=x$, ~$\partial_1 (x,y)=y$.
\item $0_x=(x,x)$,~ $1_x=(x,\theta(x))$.
\item $(x,y)^*=(y,\theta(x))$.
\item $(x,y)\cup (y,z)=(x,z)$ provided that $x\leq y \leq z\leq \theta(x)$.
\end{itemize}

Note that these definitions of $(P,S,\partial_0,\partial_1,0,1,*,\cup)$ make sense because they are compatible with the equivalence relation. Is is easy to see that one obtains in this way an abstract circle. In particular, the archimedean property \eqref{arch} is used in order to prove the Non-triviality property $1$ of an abstract circle.

Next Lemma provides the converse association.
\begin{lem}\label{equivcat0}  Let $C$ be an abstract circle then there exists an archimedean set $(X,\theta)$  such that  $C$ is isomorphic to $X/\theta$.
\end{lem}
\proof It follows from \cite{Moerdjik} that  given a point $x\in P$ the subset $L_x\subset S$ of $a\in S$ such that $\partial_0 a=x$ is an interval (with smallest element $0_x$ and largest element $1_x$) for the following order:
\begin{equation}
\label{substract}
a \leq b \iff \exists c, \  a \cup c=b.
\end{equation}
 Let then $X_x$ be the quotient of $\Z\times L_x$ by the equivalence relation which identifies $(n+1,0_x)\sim (n,1_x)$. One endows $X_x$ with the total order induced from the lexicographic ordering of $\Z\times L_x$. By construction, the order automorphism $\theta_x$ defined by $\theta_x(n,u)=(n+1,u)$ fulfills the archimedean property \eqref{arch}. The map $p: X_x\to P$, $p(n,a)=\partial_1(a)$ is a bijection of the orbit space $X_x/\sim$ (for the action of $\Z$ on $X_x$ by powers of $\theta$) with $P$. Let $s$ be the map which associates to a pair $(\alpha,\beta)\in X_x^2$ with $\alpha\leq \beta<\theta(\alpha)$, the element $s(\alpha,\beta)\in S$ determined as follows. If $\alpha=(n,a)$, $\beta=(n,b)$), one lets $s(\alpha,\beta)=c$ where $c\in S$ is determined by \eqref{substract}.  If 
$\alpha=(n,a)$, $\beta=(n+1,b)$) then  $s(\alpha,\beta)=a^*\cup b$.  Then $s$ determines a bijection of the orbit space for the diagonal action of $\Z$ on the set of pairs $(\alpha,\beta)\in X_x^2$ with $\alpha\leq \beta\leq \theta(\alpha)$ onto $S$. It is straightforward to verify that the pair $(p,s)$ provides an isomorphism $X_x/\theta_x\stackrel{\sim}{\to} C$ of abstract circles.\endproof

Let $y\in P$ be another base point, then there is a natural isomorphism $X_y\to X_x$ that is described as follows. Let $v\in S$ be the unique element such that $\partial_0 v=x$, $\partial_1 v=y$.  Let $b\in L_y$, then one has $\partial_1 v=\partial_0 b$. One sets $\psi((n,b))=(n,v\cup b)$ if $v\cup b$ exists, and 
$\psi((n,b))=(n+1,(b^*\cup v^*)^*)$ otherwise. This defines an isomorphism $\psi_{xy}:X_y\stackrel{\sim}{\to} X_x$ of archimedean sets. However it is not true that $\psi_{zy}\circ \psi_{yx}=\psi_{zx}$, in fact for $x\neq y$ one has $\psi_{xy}\circ \psi_{yx}=\theta_x$.

The above development suggests to introduce the following category
\begin{defn}\label{last} The category $\Arc$ has as objects  the archimedean sets $(X,\theta)$;  the morphisms  
$f: (X,\theta)\to (X',\theta')$ of $\Arc$ are  equivalence classes of maps 
\[
f:X\to X', \ \ f(x)\geq f(y) \  \  \forall x\geq y ;\qquad  f(\theta(x))=\theta'(f(x)),\quad \forall x \in X
\]
modulo the relation which identifies two such maps $f$ and $g$ if there exists an integer $m\in \Z$ such that 
$g(x)=\theta'^m(f(x))$, $\forall x\in X$.
\end{defn}

The full subcategory of $\Arc$ whose objects are the archimedean sets $(\Z,\theta)$, where $\Z$ is endowed with its usual order, is {\em canonically isomorphic} to the cyclic category $\Lambda$ (see Definition \ref{defncatcy}), since the archimedean automorphism $\theta$ is necessarily given by a translation $x\mapsto \theta(x)=x+n+1$, for some $n\geq 0$.
In turn, the cyclic category $\Lambda$  is isomorphic to its opposite  $\Lambda^{\rm op}$: this equivalence is established by the contravariant functor $\tt: \Lambda\longrightarrow \Lambda^{\rm op}$ which associates to $f\in \Hom_\Lambda([n],[m])$ the transposed map $f^t \in \Hom_\Lambda([m],[n])$ verifying
\[
f(x)\geq y \iff x\geq f^t(y) \qqq x,y \in \Z.
\]
The square of  $\tt$  is equivalent to the identity by means of  the natural transformation defined by the translation of $1$.
Thus to the inclusion $\Delta\subset \Lambda$ corresponds an inclusion of the opposite categories $\Delta^{\rm op}\subset \Lambda^{\rm op}\sim \Lambda$ which is described as follows: to $\varphi \in \Hom_\geq(n^*,m^*)$ one associates the unique element $\underline \varphi\in \Hom_\Lambda([n],[m])$ which satisfies $\underline \varphi(x)= \varphi(x)$, $\forall x\in \{0, \ldots, n+1\}\subset \Z$. In particular, one has: $f\in \Delta\subset \Lambda$ if and only if 
$f^t\in \Delta^{\rm op}\subset \Lambda$.

\begin{prop}\label{equivcat} 
$(i)\,$The functor $\frak Q$ which associates to an object $(X,\theta)$ of the category $\Arc$  the abstract circle $X/\theta$
is an equivalence of categories.

 $(ii)\,$The inclusion $\Delta\subset \Lambda$ induces at the level of the points of the corresponding topoi the functor which associates to an interval $I$ the archimedean set $(\Z\times I)/\sim$.
\end{prop}
\proof $(i)\,$ Lemma \ref{equivcat0} shows that any abstract circle is isomorphic to an object in the image of $\frak Q$, thus it is enough (using the definition of equivalence of categories as in \cite{MM},  p. 13) to show that $\frak Q$ is fully faithful. Let $(X,\theta)$ and $(X',\theta')$ be two archimedean sets and $h:X/\theta\to X'/\theta'$ a morphism of the associated abstract circles. Let $x\in X$ and choose $x'\in X'$ so that the image  $h(\tilde x)$ of the class $\tilde x\in P=X/\theta$ of $x\in X$ is the class $\tilde x'\in P'=X'/\theta'$.  The map $h$ determines an order preserving map from the interval $L_{\tilde x}$ (\cf \ the proof of Lemma \ref{equivcat0}) to the interval $L'_{\tilde x'}$. This map lifts to a map $f:X\to X'$ such that $\frak Q(f)=h$ and whose class in $\Hom_{\frak Q}(X,X')$ is uniquely determined. This shows that the functor $\frak Q$ is fully faithful.

 $(ii)\,$The correspondence between points of the topos of cyclic sets and abstract circles (\cite{Moerdjik}) translates, implementing $(i)$, into the following description of  points of the topos of cyclic sets. To an object $X$ of $\Arc$  one associates the restriction to $\Lambda\subset \Arc$ of the contravariant functor $\Hom_\Arc(\cdot\,, X)$ from $\Arc$  to sets. Using the contravariant functor $\tt$ one thus obtains a covariant functor $\Lambda\longrightarrow \Se$, which is also filtering. Note that all filtering functors  $\Lambda\longrightarrow \Se$  are of this form. In order to understand the effect of the inclusion $\Delta\subset \Lambda$ at the level of the points of the corresponding topoi, it is enough to consider the special case of the points associated to the objects of the (small) categories, since the other points are obtained as filtering limits of these special points. Given an object $[n]$ of $\Delta$ the associated point of the topos $\hat\Delta$ is given by the filtering (covariant) functor $\Hom_\Delta([n],\cdot\,):\Delta\longrightarrow\Se$. This is the composite of  $\Hom_\geq(\cdot\,, n^*) $ with the canonical contravariant functor $\Delta\longrightarrow\Delta^{\rm op}$. Similarly, the associated point of the topos $\hat\Lambda$ is given by the filtering (covariant) functor $\Hom_\Lambda([n],\cdot\,):\Lambda\longrightarrow\Se$ and thus it corresponds to the archimedean set $(\Z, \theta)$ where $\theta$ is the translation by $n+1$. In this dual description, the inclusion $\Delta\subset \Lambda$ corresponds  to the inclusion $\Delta^{\rm op}\subset \Lambda^{\rm op}\sim \Lambda$, where the latter is the restriction to finite intervals of the functor which associates to an interval $I$ the archimedean set $(\Z\times I)/\sim$.  \endproof


\begin{thebibliography}{99}



\bibitem{boy}  S.~Boyer,  D.~Rolfsen, B.~Wiest, {\em
Orderable 3-manifold groups}. (English, French summary) 
Ann. Inst. Fourier (Grenoble) 55 (2005), no. 1, 243--288. 

\bibitem{bu0}  D.~Burghelea, Z.~Fiedorowicz,  {\em  Cyclic homology and algebraic K-theory of spaces. II.} Topology 25 (1986), no. 3, 303--317.



\bibitem{CoExt} A.~Connes, {\em  Cohomologie cyclique et foncteurs ${\rm Ext}\sp n$},  C. R. Acad. Sci. Paris S\'er. I Math. 296 (1983), no. 23, 953--958.

\bibitem{ncg} A.~Connes, {\em Noncommutative differential geometry}, Inst. Hautes Etudes Sci. Publ.Math. No. 62 (1985), 257--360.

\bibitem{NCG} A.~Connes, {\em Noncommutative  geometry},
 Academic Press (1994).

\bibitem{cycarch} A.~Connes, C.~Consani, {\em Cyclic homology, Serre's local factors and the $\lambda$-operations}, preprint (2012).

\bibitem{FT} B.L. Feigin  and B.L. Tsygan, {\em Additive $K$-theory}.
 $K$-theory, arithmetic and geometry, pp. 67--209,  Lecture Notes
 in Math., 1289, Springer, Berlin, 1987.




\bibitem{glass} A. M. W.~Glass,  {\em Partially ordered groups}, volume 7 of Series in Algebra. World Scientific Publishing Co. Inc., River Edge, NJ, 1999.

\bibitem{good} T.~Goodwillie,   {\em Cyclic homology, derivations, and the free loopspace}. Topology 24 (1985), no. 2, 187--215. 

\bibitem{jones}  J.~Jones,   {\em Cyclic homology and equivariant homology}. 
Invent. Math. 87 (1987), no. 2, 403--423. 

\bibitem{Loday} J.L.~Loday, {\em Cyclic homology}. Grundlehren der
Mathematischen Wissenschaften, 301. Springer-Verlag, Berlin, 1998.

\bibitem{MM} S.~Mac Lane, I~Moerdijk, {\em Sheaves in geometry and logic. A first introduction to topos theory}. Corrected reprint of the 1992 edition. Universitext. Springer-Verlag, New York, 1994.

\bibitem{Moerdjik}  I~Moerdijk, {\em Cyclic sets as a classifying topos} (preprint).

\end{thebibliography}
\end{document}